\documentclass[a4paper,11pt,reqno]{amsart}
\usepackage[latin1]{inputenc}
\usepackage[english]{babel}
\usepackage{amssymb,amsfonts, amsmath, color}
\usepackage[margin=2.7cm]{geometry}
\usepackage{xcolor}
\usepackage{graphicx, color, enumerate}
\usepackage[latin1]{inputenc}
\usepackage[active]{srcltx}
\usepackage{tikz}
\usepackage{pgf}
\usepackage{etex}
\usepackage{verbatim}
\usepackage{tikz-3dplot}
\usepackage{pgfkeys}
\usepackage{amsmath}
\usepackage{amsfonts}
\usepackage{amssymb}
\usepackage{amsthm}
\usepackage{algpseudocode}
\usepackage{float}
\usepackage{xcolor}
\usepackage{mathrsfs}
\usepackage{import}
\usepackage{geometry}
\usepackage{fancyhdr}
\usepackage{fp}
\usepackage[colorlinks, citecolor=blue, linkcolor=red]{hyperref}

\newtheorem{theorem}{Theorem}[section]
\newtheorem{proposition}[theorem]{Proposition}
\newtheorem{corollary}[theorem]{Corollary}
\newtheorem{lemma}[theorem]{Lemma}
\newtheorem{remark}[theorem]{Remark}
\newtheorem{definition}[theorem]{Definition}

\newcommand{\bcl}{\begin{center}}
\newcommand{\ecl}{\end{center}}
\newcommand{\brl}{\begin{right}}
\newcommand{\erl}{\end{right}}
\newcommand{\ben}{\begin{enumerate}}
\newcommand{\een}{\end{enumerate}}
\newcommand{\overliner}{\begin{array}}
\newcommand{\earr}{\end{array}}
\newcommand{\btab}{\begin{tabular}}
\newcommand{\etab}{\end{tabular}}
\newcommand{\bdoc}{\begin{document}}
\newcommand{\edoc}{\end{document}}
\newcommand{\beqy}{\begin{eqnarray}}
\newcommand{\eeqy}{\end{eqnarray}}

\newcommand{\beqi}{\begin{eqnarray*}}
\newcommand{\eeqi}{\end{eqnarray*}}
\newcommand{\bitem}{\begin{itemize}}

\newcommand{\eitem}{\end{itemize}}
\newcommand{\nln}{\newline}
\newcommand{\newt}{\newtheorem}

\newcommand{\pa}{\partial}
\newcommand{\re}{{I\!\!R}}
\newcommand{\Rn}{\R^N}
\newcommand{\xr}{x\in\R }
\newcommand{\x}{\times}
\newcommand{\dyle}{\displaystyle}
\newcommand{\ene}{{I\!\!N}}
\newcommand{\irn}{\int\limits_{\R^N}}
\newcommand{\io}{\int\limits_{\O}}
\newcommand{\meas}{{\rm meas\,}}
\newcommand{\dif}{\nabla_{xy}}
\newcommand{\sign}{{\rm sign}}
\newcommand{\map}{\longrightarrow }
\newcommand{\imp}{\Longrightarrow }
\renewcommand{\div}{\nabla\cdot }
\newcommand{\sen}{{\rm sen\,}}
\newcommand{\tg}{{\rm tg\,}}
\newcommand{\arcsen}{{\rm arcsen\,}}
\newcommand{\arctg}{{\rm arctg\,}}
\newcommand{\supp}{{\textsl supp\ }}
\newcommand{\ity}{\int_{-\iy}^{+\iy}}
\newcommand{\limit}{\lim\limits}
\newcommand{\limi}{\limit_{n\to\infty}}
\newcommand{\sumi}{\sum\limits_{n=1}^{\infty}}
\newcommand{\ulu}{\underline u}
\newcommand{\ulw}{\underline w}
\newcommand{\ulz}{\underline z}
\newcommand{\ulv}{\underline v}
\newcommand{\uls}{\underline s}
\newcommand{\olu}{\overline u}
\newcommand{\olv}{\overline v}
\newcommand{\ols}{\overline s}
\newcommand{\ob}{\overline\b}
\newcommand{\ovar}{\overline\var}
\newcommand{\wv}{\widetilde v}
\newcommand{\wu}{\widetilde u}
\newcommand{\ws}{\widetilde s}
\renewcommand{\a }{\alpha }
\renewcommand{\b }{\beta }
\newcommand{\g }{\gamma}
\newcommand{\G }{\Gamma }
\renewcommand{\d }{\delta }

\newcommand{\D }{\Delta }
\newcommand{\e }{\varepsilon }
\newcommand{\z }{\zeta }
\renewcommand{\l }{\lambda }
\renewcommand{\L }{\Lambda }
\newcommand{\m }{\mu }
\newcommand{\n }{\nabla }
\newcommand{\s }{\sigma }
\newcommand{\Sig }{\Sigma }
\renewcommand{\t }{\tau }
\newcommand{\var }{\varphi }
\renewcommand{\o }{\omega }
\renewcommand{\O }{\Omega }
\newcommand{\R}{{\mathbb{R}}}
\newcommand{\bC}{{\bf C}}
\newcommand{\bZ}{{\bf Z}}
\newcommand{\bN}{{\bf N}}
\newcommand{\bQ}{{\bf Q}}
\newcommand{\bK}{{\bf K}}
\newcommand{\bI}{{\bf I}}
\newcommand{\bv}{{\bf v}}
\newcommand{\bV}{{\bf V}}
\newcommand{\LL}{\mathcal{L}}
\newcommand{\N}{\mathbb{N}}
\DeclareMathOperator{\suppo}{supp} \DeclareMathOperator{\di}{div}

\newenvironment{Proof}{\Rmovelastskip\vskip12pt
plus 1pt \noindent\em\rm}{\hfill {\qed \hskip .2cm}}

\begin{document}

\title[]{A Liouville theorem \\ for  elliptic equations with a potential \\ on infinite graphs}

\author{Stefano Biagi}

\address{\hbox{\parbox{5.7in}{\medskip \noindent{Stefano Biagi, \\Dipartimento di Matematica, \\Politecnico di Milano, \\Piazza Leonardo da Vinci 32, 20133, Milano, Italy \\ [3pt] \emph{E-mail address: }{\tt stefano.biagi@polimi.it}}}}}

\author{Giulia Meglioli}

\address{\hbox{\parbox{5.7in}{\medskip \noindent{Giulia Meglioli, \\Fakult\"at f\"ur Mathematik, \\Universit\"at Bielefeld, \\33501, Bielefeld, Germany \\ [3pt] \emph{E-mail address: }{\tt gmeglioli@math.uni-bielefeld.de}}}}}

\author{Fabio Punzo}

\address{\hbox{\parbox{5.7in}{\medskip \noindent{Fabio Punzo, \\Dipartimento di Matematica, \\Politecnico di Milano, \\Piazza Leonardo da Vinci 32, 20133, Milano, Italy \\ [3pt] \emph{E-mail address: }{\tt fabio.punzo@polimi.it}}}}}



\begin{abstract} We investigate the validity of the Liouville property for a class of elliptic equations with a potential, posed on infinite graphs. Under suitable assumptions on the graph and on the potential, we prove that the unique bounded solution is $u\equiv 0$. We also show that on a special class of graphs the condition on the potential is optimal, in the sense that if it fails, then there exist infinitely many bounded solutions.
\end{abstract}

\maketitle

\section{Introduction} \label{sec0}

Let
$(G, \omega, \mu)$ be
 a fixed {\it infinite weighted} graph, with {\it edge-weight} $\omega$ and {\it node (or vertex) measure} $\mu$.
In this paper we study {\it bounded} solutions to elliptic equations with a potential of the following form:
\begin{equation}\label{problema}
\Delta u -V u=0 \quad \text{ in }\; G,
\end{equation}
where the potential $V$ is a nonnegative function defined in $G$
and $\Delta$ denotes the Laplace operator on $G$.

The uniqueness of solutions of equation \eqref{problema} has been
investigated by two of us in the recent paper \cite{MP2}; in this paper,
it is proved that $u\equiv 0$ is the only solution  to equation \eqref{problema}, whenever that $u$ belongs to a certain $\ell^p_{\varphi}(G, \mu)$ space, where $\varphi$ is a weight which tends to $0$ at {\it infinity} and $p\in [1, +\infty)$. Also the case $u\in \ell^\infty(G, \mu)$ can be considered, provided that the graph satisfies a suitable property.  In any case, an essential hypothesis
for the arguments used in \cite{MP2} is the existence of some $c_1 > 0$ such that
\begin{equation}\label{e301f}
V(x)\geq c_1 \quad \text{ for all }\; x\in G\,.
\end{equation}
Hence, the aim of the present paper is to understand when $u\equiv 0$ is the  unique {\it bounded} solution of problem \eqref{problema}, without supposing hypothesis \eqref{e301f}. To do this, we have to use completely different methods than those exploited in \cite{MP2}.

\smallskip

We say that the {\it Liouville theorem (or property)} holds for equation \eqref{problema}, whenever $u\equiv 0$ is the only {\it bounded} solution of the same equation. Thus, in other terms, we are concerned with the validity of the Liouville theorem for equation \eqref{problema}.
\smallskip

Before describing our results and methods of proof, let us contextualize our problem within the literature. As it is well-known, many phenomena in various fields of applied sciences can be modeled by means of graphs (see, e.g., \cite{Deo, Les, LHN})\,. For this reason, partial differential equations posed on graphs have
recently attracted the attention of many authors. In particular, qualitative properties of solutions  to both elliptic and parabolic equations have been addressed, see {\it e.g.}  \cite{GLY1, GLY2, GLY3, GHY, HW, KS, LY, PinaS} and  \cite{BCG, CGZ, EM, HL, HMu, LW2, Mu, SSV, Wu}, respectively. Moreover, the monographs \cite{Grig2, KLW, Mu2} are important contributions to this topic.

In \cite{Huang} and in  \cite{HKS}, under suitable assumptions on the graph, it is shown that the parabolic problem
\[\begin{cases}
\partial_t u = \Delta u  & \text{ in } G\times (0, \infty)\\
u = 0 & \text{ in } G\times \{0\}\,
\end{cases} \]
has at most one solution fulfilling a suitable growth condition. An analogous result can be found in \cite{CGZ}, where the time derivative is replaced by discrete differences. Furthermore, in \cite[Theorem 12.15, Corollary 12.16, Theorem 12.17]{KLW} it is shown that, under suitable assumptions on $G$, if $u$ is a
subsolution
of equation \eqref{problema} with $V\equiv 0$ which satisfies
 $u\in \ell^p(G, \mu)$ and $u\geq 0$, then $u$ must be constant.

Similar uniqueness results have been established also on manifolds (see, e.g., \cite{Grig3, Grig, GrigHK, MR, Pu1}), in bounded domains of $\mathbb R^n$ (see \cite{BP1, BP2, NP, Pu2}), and for nonlocal operators (see \cite{MP, PV}).

\medskip

Now, let $d$ denote a {\it pseudo metric} on $G$, and let $B_r(x_0)$ be the ball centered at $x_0\in G$ with radius $r$ (see Section \ref{mf} below). Concerning the potential $V$, we suppose that
\[V(x)\geq c_0\, d^{-\alpha}(x, x_0) \quad \textrm{for all}\;\; x\in G\setminus B_{R_0}(x_0),\]
for some $c_0>0, R_0>0$ and $\alpha\in [0, 1]$.
Under this assumption on $V$ we prove that, if there exists some
constant $\Lambda\in (0,1)$ such that
 \begin{equation}\label{eqag1}
  \sum_{x\in G\setminus B_{1}(x_0)}e^{-\Lambda d^\alpha(x,x_0)}\mu(x) < +\infty,
  \end{equation}
then $u\equiv 0$ is the only bounded solution of equation \eqref{problema}.

In order to prove this result, we take inspiration from \cite{Grig3}. However, since in \cite{Grig3} the problem is posed on a Riemannian manifold, many ideas used in that setting cannot be exploited on graphs, hence important differences arise.
More precisely, the line of arguments we follow to show our main result is the following:
introducing the function
\[v(x,t):=e^{t}u(x) - 1 \quad \text{for all }\; x\in G,\,\, t\in [0, T],\]
we first show that the positive part of $v$, namely $v_+$, is a subsolution of equation \eqref{problema} in an appropriate sense (see Lemma \ref{lem4}). Then, we obtain a key {\it a priori} estimate for $v_+$ (see Proposition \ref{lemma1}), where test functions $\xi=\xi(x,t)$ and $\eta=\eta(x)$ are employed.
Then we select a suitable function $\xi$, which can be regarded as a sort of supersolution of an ``adjoint parabolic equation" (see Lemma \ref{lemma2}). On the other hand, $\eta$ will be chosen to be a ``cut-off" function (see Lemma \ref{lemma3}).
Then the conclusion follows by means of appropriate estimates.

A key point in our strategy of proof of uniqueness is to show that if $u$ is a bounded solution of equation \eqref{problema}, then there exists a solution $v$ to the same problem such that $0<v\leq 1$ (see Proposition \ref{prop2}).
To show this, we first need to establish
a weak maximum principle, see Lemma \ref{lem:WMP}, and a strong maximum principle, see Lemma \ref{lem:SMP}. Let us mention that Proposition \ref{prop2}, Lemmas  \ref{lem:WMP}, \ref{lem:SMP}, \ref{lem4} can have an independent interest.

\smallskip

We also show that the bound \eqref{eqag1} with $\alpha\in [0, 1]$ is optimal. More precisely, on a special class of graphs, we prove that if $V$ decays like $d^{-\alpha}(x, x_0)$ for some $\alpha>1$ as $d(x, x_0)\to +\infty$, then there exist infinitely many bounded solutions to problem \eqref{problema}. To be specific, for any $\gamma>0$ there exists a $u$ solution to \eqref{problema} such that $u(x)\to \gamma$ as $d(x, x_0)\to +\infty$. Its proof is based on the construction of a suitable barrier at infinity, which is  related to the class of graphs we consider.
To the best of our knowledge, on graphs, such type of result, which consists in prescribing a Dirichlet condition at {\it infinity}, and the explicit construction of such kind of barrier are totally new.
\medskip

The paper is organized as follows. In Section \ref{mf} we describe the graph framework and we state our main result. Section \ref{auxiliary} is devoted to the auxiliary results for elliptic equations previously described. In Section \ref{apriori} we obtain the key apriori estimate for $v_+(x,t)$. We introduce and study our test functions in Section \ref{test}. The main result is proved in Section \ref{proofthm}. Finally, in Section \ref{optimality} we show the optimality of condition \eqref{eqag1}.

\section{Mathematical framework and the main result}\label{mf}\setcounter{equation}{0}

\subsection{The graph setting}
Let $G$ be a countably infinite set and $\mu:G\to (0,+\infty)$ be a given function. Observe that $\mu$ can be viewed as a Radon measure on $G$ so that $(G,\mu)$ becomes a measure space. Furthermore, let
\begin{equation*}
\omega:G\times G\to [0,+\infty)
\end{equation*}
be a symmetric, with zero diagonal and finite sum function, i.e.
\begin{equation}\label{omega}
\begin{aligned}
&\text{(i)}\,\, \omega(x,y)=\omega(y,x)\quad &\text{for all}\,\,\, (x,y)\in G\times G;\\
&\text{(ii)}\,\, \omega(x,x)=0 \quad\quad\quad\,\, &\text{for all}\,\,\, x\in G;\\
&\text{(iii)}\,\, \displaystyle \sum_{y\in G} \omega(x,y)<\infty \quad &\text{for all}\,\,\, x\in G\,.
\end{aligned}
\end{equation}
Thus, we define  \textit{weighted graph} the triplet $(G,\omega,\mu)$, where $\omega$ and $\mu$ are the so called \textit{edge weight} and \textit{node measure}, respectively. Observe that assumption $(ii)$ corresponds to ask that $G$ has no loops.
\smallskip

\noindent Let $x,y$ be two points in $G$; we say that
\begin{itemize}
\item $x$ is {\it connected} to $y$ and we write $x\sim y$, whenever $\omega(x,y)>0$;
\item the couple $(x,y)$ is an {\it edge} of the graph and the vertices $x,y$ are called the {\it endpoints} of the edge whenever $x\sim y$;
\item a collection of vertices $ \{x_k\}_{k=0}^n\subset G$ is a {\it path} if $x_k\sim x_{k+1}$ for all $k=0, \ldots, n-1.$
\end{itemize}

\noindent We are now ready to list some properties that the weighted graph $(G,\omega,\mu)$ may satisfy.
\begin{definition}\label{def01}
We say that the weighted graph $(G,\omega,\mu)$ is
\begin{itemize}
\item[(i)] {\em locally finite} if each vertex $x\in G$ has only finitely many $y\in G$ such that $x\sim y$;
\item[(ii)] {\em connected} if, for any two distinct vertices $x,y\in G$ there exists a path joining $x$ to $y$;
\item[(iii)] {\em undirected} if its edges do not have an orientation.
\end{itemize}
\end{definition}
For any $x\in G$, we define
\begin{itemize}
\item the {\it degree} of $x$ as $$\operatorname{deg}(x):=\sum_{y\in G}\omega(x,y);$$
\item the {\it weighted degree} of $x$ as $$\operatorname{Deg}(x):=\frac{\operatorname{deg}(x)}{\mu(x)}.$$
\end{itemize}

\noindent A {\it pseudo metric} on $G$ is a symmetric, with zero diagonal map, $d:G\times G\to [0, +\infty)$, which also satisfies the triangle inequality
\begin{equation*}
  d(x,y)\leq d(x,z)+d(z,y)\quad \text{for all}\,\,\, x,y,z\in G.
\end{equation*}
In general, $d$ is not a metric, since we can find points $x, y\in G$, $x\neq y$ such that $d(x,y)=0\,.$ Now, let us consider any path $\gamma\equiv \{x_k\}_{k=0}^n\subset G$, and a symmetric function $\sigma:G\times G\to [0, +\infty)$ such that
$$
\sigma(x,y)>0 \quad \text{if and only if}\,\,\, x\sim y.
$$
Then we define the {\it lenght subordinated to} $\sigma$ as
$$
l_\sigma(\gamma):=\sum_{k=0}^{n-1}\sigma(x_k, x_{k+1})\,.
$$
A pseudo metric $d\equiv d_\sigma$ is called a {\it path pseudo metric} if there exists a symmetric map $\sigma:G\times G\to [0, +\infty)$, with $\sigma(x,y)>0$ if and only if $x\sim y$, such that
\begin{equation*}
d_{\sigma}(x,y)=\inf\left\{\,l_{\sigma}(\gamma) : \gamma\,\, \text{is a path between}\,\, x\,\, \text{and}\,\, y\,\right\}.
\end{equation*}
Finally, we define the \textit{jump size} $s>0$ of a pseudo metric $d$ as
\begin{equation}\label{e14f}
s:=\sup\{d(x,y) \,:\, x,y\in G, \omega(x,y)>0\}.
\end{equation}
For a more detailed understanding of the objects introduced so far, we refer the reader to \cite{GT,HKMW,HKS,MP2}. We conclude the subsection with the following
\begin{definition}\label{intrmet}
A pseudo metric $d$ on $(G, \omega, \mu)$ is said to be {\em intrinsic} if
\begin{equation}\label{e1f}
\frac 1{\mu(x)}\sum_{y\in G} \omega(x,y)d^2(x,y) \leq 1\quad \text{ for all } x\in G\,.
\end{equation}
\end{definition}
\noindent Observe that hypothesis \eqref{e1f} can be compared with an analogous
condition on Riemannian manifolds. Indeed,
given any fixed reference point $x_0\in G$, let us
consider the map
$$G\ni x\mapsto d(x, x_0).$$
Then, we have
\begin{align*}
 |\nabla d(x, x_0)|^2 & =\frac 1{\mu(x)}\sum_{y\in G}\omega(x,y)[d(y, x_0)-d(x, x_0)]^2 \\
 & \leq \frac 1{\mu(x)} \sum_{y\in G}\omega(x,y) d^2(x,y)\;  \text{ for any } x\in G.
 \end{align*}
Therefore, condition \eqref{e1f} ensures that
\[|\nabla d(x, x_0)|^2\leq 1 \quad \text{ for any } x\in G\,.\]
Such a property is clearly fulfilled on Riemannian manifolds.
\medskip

For any $x_0\in G$ and $r>0$ we define the ball $B_r(x_0)$ with respect to any pseudo metric $d$ as
\[B_r(x_0):=\{x\in G\,:\, d(x,x_0)<r\}\,.\]

\noindent In this paper, we always make the following assumptions:
\begin{equation}\label{e7f}
\begin{aligned}
(i)\,\,\, & (G, \omega, \mu) \text{ is a connected, locally finite, weighted graph};\\
(ii) \,\,\,& \text{there exists a \textit{pseudo metric}}\,\, d \,\,\,\text{such that the jump size $s$ is finite}; \\
(iii)\,\,\,& \text{the ball}\,\,\, B_r(x) \,\,\,\text{with respect to}\,\,\, d\,\,\, \text{is a finite set, for any}\,\,\, x\in G,\,\,\, r>0;\\
\end{aligned}
\end{equation}
here we have used Definitions \ref{def01} and \ref{intrmet}.

\subsection{Difference and Laplace operators} Let $\mathfrak F$ denote the set of all functions $f: G\to \mathbb R$\,. For any $f\in \mathfrak F$ and for all $x,y\in G$, let us give the following
\begin{definition}\label{def1}
Let $(G, \omega,\mu)$ be a weighted graph. For any $f\in \mathfrak F$,
\begin{itemize}
\item the {\em difference operator} is
\begin{equation}\label{e2f}
\nabla_{xy} f:= f(y)-f(x)\,;
\end{equation}
\item the {\em (weighted) Laplace operator} on $(G, \omega, \mu)$ is
\begin{equation*}
\Delta f(x):=\frac{1}{\mu(x)}\sum_{y\in G}[f(y)-f(x)]\omega(x,y)\quad \text{ for all }\, x\in G\,.
\end{equation*}
\end{itemize}
\end{definition}
Clearly,
\[\Delta f(x)=\frac 1{\mu(x)}\sum_{y\in G}(\dif f)\omega(x,y)\quad \text{ for all } x\in G\,.\]
We also define  the {\it gradient squared} of $f\in \mathfrak F$ (see \cite{CGZ})
\[|\nabla f(x)|^2=\frac 1{\mu(x)}\sum_{y\in G}\omega(x,y)(\dif f)^2, \quad x\in G\,;\]
It is straightforward to show, for any $f,g\in \mathfrak F$, the validity of
\begin{itemize}
\item  the {\it product rule}
\begin{equation*}
\nabla_{xy}(fg)=f(x) (\nabla_{xy} g) + (\nabla_{xy} f)g(y) \quad \text{ for all } x,y\in G\,;
\end{equation*}
\item the {\it integration by parts formula}
\begin{equation}\label{e4f}
\sum_{x\in G}[\Delta f(x)] g(x) \mu(x)=-\frac 1 2\sum_{x,y\in G}(\dif f)(\dif g) \omega(x,y)\,,
\end{equation}
provided that at least one of the functions $f, g\in \mathfrak F$ has {\it finite} support.
\end{itemize}



\subsection{The main result} \label{sec2}

\noindent We have already stated in \eqref{e7f}
the main hypotheses on the weighted graph $(G,\omega,\mu)$.
Concerning the {\it potential} $V$, we  suppose that
\begin{equation}\label{e12f}
\begin{aligned}
&V\in \mathfrak F,\\
& V(x)\geq 0 \quad \text{ for all } x\in G,\\
& V(x)\geq c_0\, d^{-\alpha}(x, x_0) \quad \textrm{for all}\;\; x\in G\setminus B_{R_0}(x_0)\,.
\end{aligned}
\end{equation}
for some $x_0\in G, R_0>0$, $c_0>0$ and $\alpha\in [0, 1].$
\vspace*{0.1cm}

We can now state the main result of this paper.

\begin{theorem}\label{teo1}
Let assumptions \eqref{e1f}, \eqref{e7f} and \eqref{e12f} be satisfied.  Let $u$ be a bounded solution of problem \eqref{problema}. Moreover, suppose that there exists $\Lambda \in (0,1)$ such that
\begin{equation} \label{eq:hpStrada3}
  \sum_{x\in G\setminus B_{1}(x_0)}e^{-\Lambda d^\alpha(x,x_0)}\mu(x) < +\infty.
  \end{equation}
Then
$$u(x) = 0 \quad \text{ for any } \; x\in G.$$
\end{theorem}

\begin{remark} Observe that
condition \eqref{eqag1} with $\alpha\in [0, 1]$ is sharp. In fact, on a special class of graphs, the so-called model trees, fulfilling \eqref{eq:hpStrada3}, if $V$ decays like $d^{-\alpha}(x, x_0)$ for some $\alpha>1$ as $d(x, x_0)\to +\infty$, then there exist infinitely many bounded solutions to problem \eqref{problema}. More precisely, for any $\gamma\in (0, +\infty)$ there exists a solution $u$ to \eqref{problema} such that 
\[u(x)\to \gamma \quad \text{as }\; d(x,x_0)\to +\infty\,.\]
Since $\gamma>0$ is arbitrary, in particular nonuniqueness for equation \eqref{problema} follows. See Section \ref{optimality} for more details. To the best of our knowledge, in the literature such methods used on graphs cannot be found.
\end{remark}

\section{Auxiliary Results}\label{auxiliary}\setcounter{equation}{0}
In this section we collect several preliminary results of independent interest which shall be used in the proof of our main result, namely Theorem \ref{teo1}.
\begin{proposition}\label{prop2} Let assumption \eqref{e7f}-(i) be fulfilled. Assume that there exists a nontrivial bounded solution of equation \eqref{problema}. Then there exists a solution $v$ of equation \eqref{problema} such that
\begin{equation}\label{e1g}
0<v \leq 1 \quad \text{ in }\; G\,.
\end{equation}
\end{proposition}
Analogously to \cite{BP1, Grig3},
the proof of Proposition \ref{prop2} is crucially based on the \emph{unique solvability}
of the Dirichlet problem for $$\mathcal{L} := \Delta-V(x),$$ that is,
\begin{equation}\label{D}
\begin{cases}
\LL u = f & \text{in $\Omega$} \\
u \equiv g & \text{in $G\setminus\Omega$}
\end{cases}
\end{equation}
(where $\Omega\subseteq G$ is an arbitrary finite set and $f,g\in\mathfrak{F}$),
together with some \emph{maximum principles} for $\mathcal{L}$.
Since we were not able to find a precise reference for these results,
and in order to make the paper as self-contained as possible, we present here below
the full proofs.
\vspace*{0.1cm}

To begin with, we give the following definition
\begin{definition}\label{defsol}
We say that $u\in\mathfrak{F}$ is a solution of equation \eqref{D} if,
\begin{equation}\label{sol}
\frac{1}{\mu(x)}\sum_{y\in \Omega}\omega(x,y)\left[u(y)-u(x)\right]-V(x)u(x)=f(x)\quad\text{for every}\,\,x\in\Omega,
\end{equation}
and $u \equiv g$ in $G\setminus\Omega$. Moreover, we say that $u$ is a supersolution (subsolution) to equation \eqref{D}, if the $=$ in \eqref{sol} is replaced by $\le$ ($\ge$) and $u\ge g$ ($\le$) in $G\setminus\Omega$.
\end{definition}
We now establish the following \emph{Weak Maximum Principle}.
\begin{lemma} \label{lem:WMP} Let assumption \eqref{e7f}-(i) be fulfilled.
 Let $\Omega\subseteq G$ be a finite set, and let $u\in\mathfrak{F}$ be such that
 \begin{equation} \label{eq:PBWMP}
  \begin{cases}
   \LL u \leq 0 & \text{in $\Omega$} \\
   u \geq 0 & \text{in $G\setminus\Omega$}.
  \end{cases}
 \end{equation}
 Then $$u\geq 0 \quad \text{ in } \Omega.$$
\end{lemma}
\begin{proof}
 We proceed essentially as in the proof of \cite[Lemma 1.39]{Grig2}. We set $m := \min_\Omega u$; observe that $m$ is well-defined since the set $\Omega\subseteq G$ is finite.
Suppose, by contradiction, that $m < 0$. Then the set
 $$F := \{x\in G:\,u(x) = m\}\neq\varnothing$$
 is such that
 \begin{equation} \label{eq:connectedSWMP}
  \text{if $x\in F$ and $y\in G,\,y\sim x$, then $y\in F$}.
 \end{equation}
Indeed, let $x\in F$ be fixed, hence $u(x)=m < 0$. Due to \eqref{e12f}, \eqref{eq:PBWMP} and recalling that $\omega(x,y) > 0$ if $y\sim x$, we have
 \begin{align*}
  0 \geq \LL u(x) &= \Delta u(x)-V(x)u(x) \\
  &= \frac{1}{\mu(x)}\sum_{y\in G}\omega(x,y)(u(y)-u(x))-V(x)u(x) \\
  & = -\mathrm{Deg}(x)u(x)+\frac{1}{\mu(x)}\sum_{y\sim x}\omega(x,y)u(y) -V(x)u(x) \\
  & \geq -\mathrm{Deg}(x)u(x)+\frac{1}{\mu(x)}\sum_{y\sim x}\omega(x,y)u(y).
 \end{align*}
 Therefore, since $u\geq m$ in $\Omega$ and $u\geq 0 > m$ in $G\setminus \Omega$, we obtain
 $$m\mathrm{Deg}(x) = \mathrm{Deg}(x)u(x) \geq \frac{1}{\mu(x)}
   \sum_{y\sim x}\omega(x,y)u(y) \geq
  m\mathrm{Deg}(x),$$
 from which we derive that
 \begin{equation} \label{eq:ineqeqWMP}
  \sum_{y\sim x}\omega(x,y)u(y) = m.
  \end{equation}
  In view of \eqref{eq:ineqeqWMP}, since $u\geq m$ in $G$, we conclude that $u(y) = m$ for every $y\in G,\,y\sim x$, i.e. \eqref{eq:connectedSWMP}.
 \smallskip

 Now, let us consider some $x\in F$ and $y\in G\setminus \Omega$, hence $u(x)=m<0$ and $u(y)\ge 0$. Due to \eqref{e7f}, there exist a path $\{x_k\}_{k=0}^n$ such that
 $$x_0 = x,\,\,\,x_n = y.$$
 Since $x_0 = x\in F$, we can apply \eqref{eq:connectedSWMP} and infer that $x_1\in F$. By repeating this argument, we get that $x_i\in F$ for every $i=0,...,n$, hence in particular that $x_n = y\in F$ and thus $u(y)=m<0$ which yields a contradiction.
 \end{proof}

 We now prove the following \emph{Strong Maximum Principle} for $\LL$-harmonic functions.
 \begin{lemma} \label{lem:SMP} Let assumption \eqref{e7f}-(i) be fulfilled.
  Let $u\in\mathfrak{F}$ be such that
  \begin{equation}\label{eq381}
  \begin{cases}
  \LL u \leq 0 & \text{in $G$}, \\
  u \geq 0 & \text{in $G$}.
  \end{cases}
  \end{equation}
  Then either $u\equiv 0$ in $G$ or $u > 0$ in $G$.
 \end{lemma}
 \begin{proof}
Let $u$ be a solution of \eqref{eq381}. Then, due to Lemma \ref{lem:WMP}, $u\ge0$ in $G$. Let us now assume that there exists $x_0\in G$ such that $u(x_0) = 0$. Moreover, we consider the set
$$
F := \{x\in G:\,u(x) = 0\}\subseteq G.
$$
Observe that $x_0\in F$. By arguing as in the proof of Lemma \ref{lem:WMP}, by using that $\min_G u=0$, we get, equivalently to \eqref{eq:connectedSWMP}, that
$$
\text{if}\,\,\, y\in G\,\,\, \text{and}\,\,\, y\sim x_0,\,\,\, \text{then}\,\,\, y\in F.
$$
Consequently, since $G$ is connected, we conclude that $F = G$, and hence $u\equiv 0$ in $G$.
\end{proof}
Due to Lemmas \ref{lem:WMP} and \ref{lem:SMP}, we can now prove the following
\begin{proposition} \label{prop:SolvDir} Let assumption \eqref{e7f}-(i) be fulfilled.
 Let $\Omega\subseteq G$ be a finite set. Let $f:\Omega\to\R$ and $g:G\setminus\Omega\to\R$
 be \emph{arbitrary} functions. Then there exists a \emph{unique solution} $u\in\mathfrak{F}$ to problem \eqref{D} in the sense of Definition \ref{defsol}.
\end{proposition}
\begin{proof} We begin by observing that, given any $u\in\mathfrak{F}$, we can write
\begin{equation} \label{eq:rewriteLLu}
\begin{split}
 \LL u(x) & = \Delta u(x)-V(x)u(x) = \frac{1}{\mu(x)}\sum_{y\in G}\omega(x,y)(u(y)-u(x))-V(x)u(x) \\
 & = -\mathrm{Deg}(x)u(x)+\sum_{y\in G}\frac{\omega(x,y)}{\mu(x)}u(y)-V(x)u(x) \\
 & = -(\mathrm{Deg}(x)+V(x))u(x)+\sum_{y\in \Omega}P(x,y)u(y)+\sum_{y\notin\Omega}P(x,y)u(y),
\end{split}
\end{equation}
where
$$P(x,y) := \frac{\omega(x,y)}{\mu(x)}.$$
Thus, we see that $u\in\mathfrak{F}$ is a solution of \eqref{D} \emph{if and only if}
$u\equiv g$ in $G\setminus \Omega$ and
\begin{equation} \label{eq:identityAbstract}
 -(\mathrm{Deg}(x)+V(x))u(x)+\sum_{y\in \Omega}P(x,y)u(y) = f(x)-\sum_{y\notin\Omega}P(x,y)g(y)\quad
\forall\,\,x\in\Omega.
\end{equation}
 We now claim that the validity of \eqref{eq:identityAbstract}, which only involves the
 values attained by $u$ on $\Omega$, can be rephrased as a
 linear equation
 in a suitable finite-dimensional vector space.

 In fact, if we denote by $\mathfrak{F}_\Omega$ the set of all real-valued functions defined on $\Omega$,
 it is immediate to recognize
 that $\mathfrak{F}_\Omega$ is a real vector space, and
 $$\mathfrak{F}_\Omega = \mathrm{span}\{\chi_{\{x\}}:\,x\in\Omega\}$$
 (here, $\chi_A$ stands for the indicator function of the set $A\subseteq G$). Thus, since
 $\Omega$ is finite, we derive that
 $\mathfrak{F}_\Omega$ has finite dimension $n = \mathrm{card}(\Omega)$.
 On this space $\mathfrak{F}_\Omega$, we then define the map
 $$\mathcal{A}:\mathfrak{F}_\Omega\to \mathfrak{F}_\Omega,\qquad
 (\mathcal{A}u)(x) := -(\mathrm{Deg}(x)+V(x))u(x)+\sum_{y\in \Omega}P(x,y)u(y).$$
 Clearly, $\mathcal{A}$ is linear; moreover, identity \eqref{eq:identityAbstract} can be rewritten as
 $$\mathcal{A}u = h_{f,g},
  \qquad\text{where $h_{f,g} = f(x)-\sum_{y\notin\Omega}P(x,y)g(y)\in\mathfrak{F}_\Omega$}.$$
 Summing up, we have that $u\in\mathfrak{F}$ is a solution of problem \eqref{D} \emph{if and only if}
 \begin{equation} \label{eq:abstractDP}
  \text{$u\equiv g$ in $G\setminus\Omega$}\qquad\text{and}\qquad
 \mathcal{A}(u|_\Omega) = h_{f,g}.
 \end{equation}
 Using this `abstract' formulation of the Dirichlet problem \eqref{D},
 we can easily complete the proof of the proposition.
 First of all we observe that, owing to the Weak Maximum Principle
 in Lemma \ref{lem:WMP}, the linear operator $\mathcal{A}$ is injective: indeed,
 if $u\in \mathfrak{F}_\Omega$ is such that $\mathcal{A}u = 0$
 (that is, $\mathcal{A}u(x) = 0$ for every $x\in\Omega$)
 and if define $\hat{u} := u\mathbf{1}_{\Omega}$, from \eqref{eq:rewriteLLu} we have
 $$\text{$\LL \hat{u} = \mathcal{A}u = 0$ in $\Omega$};$$
 thus, since $\hat{u}\equiv 0$ in $G\setminus\Omega$, an immediate
 application of Lemma \ref{lem:WMP}
 shows that $\hat{u}\equiv 0$ in $G$, and hence $u = 0$ in $\mathfrak{F}_\Omega$.
 From this, since $\mathfrak{F}_\Omega$ has finite dimension, we derive that
 $\mathcal{A}$ is also surjective, and thus there exists a \emph{unique}
 function $u\in\mathfrak{F}_\Omega$ such that
 \begin{equation} \label{eq:toconcludeAh}
  \mathcal{A}u = h_{f,g}.
 \end{equation}
 Extending this unique function $u$ by setting $u(x) = g(x)$ for every $x\in G\setminus\Omega$,
 from \eqref{eq:toconcludeAh} we con\-clude
 that \eqref{eq:abstractDP} is satisfied, and thus $u$ is the unique solution of \eqref{D}.
\end{proof}
 With Proposition \ref{prop:SolvDir} at hand, we can finally prove Proposition \ref{prop2}.
 \begin{proof}[Proof of Proposition \ref{prop2}]
  Let $u\in\mathfrak{F}$ be a non-trivial bounded solution of problem \eqref{problema},
  and let
  $\mathbf{s} := \sup_{G}|u|\in (0, +\infty)$. Without lost of generality, we may assume $\mathbf{s} = 1$, indeed, it would be sufficient to replace $u$ with $\frac{u}{\mathbf{s}}$. Moreover, if $u$ has \emph{constant sign} on $G$, then the function
  $$v := \mathrm{sgn}(u)u,$$
  is a solution of \eqref{problema} satisfying \eqref{e1g}.
  Indeed, it is immediate to recognize that
  $v$ is a non-trivial solution of \eqref{problema}, and $0\leq v\leq 1$ on $G$. Then,
  by the Strong
  Maximum Principle in Lemma \ref{lem:SMP}, we conclude that $v > 0$ on $G$.
  If, instead, $u$ \emph{changes sign} in $G$,
  we fix $M$ such that
\[0<M<\sup_{G} u_+ = 1,\]
 and, for every  $n\in \mathbb N$, we let $v_n\in\mathfrak{F}$ be the unique solution of problem
\begin{equation*}
\begin{cases}
\LL v = 0 & \text{in $B_n(p)$}, \\
v = (u-M)_+ & \text{in $G\setminus B_n(p)$},
\end{cases}
\end{equation*}
where $p\in G$ is a point arbitrarily chosen.
The existence and uniqueness of $v_n$ for each $n\in \mathbb N$ is guaranteed by Proposition
\ref{prop:SolvDir}, since the balls $B_r(x)$ are \emph{finite sets}, see
\eqref{e7f}.

We now claim that, for every $n\in\N$, the following properties holds:
\begin{itemize}
 \item[(i)] $0\leq v_n\leq 1$ pointwise on $G$;
 \item[(ii)] $v_n\leq v_{n+1}$ pointwise on $G$.
\end{itemize}
Taking this claim for granted for a moment, we can easily complete the proof of the proposition.
In fact, owing to (i)-(ii) we deduce that the sequence $\{v_n\}_n$ is \emph{increasing and bounded}
on $G$; as a consequence, the function
$$v(x) = \lim_{n\to+\infty}v_n(x)\qquad x\in G,$$
is well-defined, and it satisfies $0\leq v\leq 1$ on $G$. Moreover, since
\begin{equation} \label{eq:LLvzero}
\begin{split}
 \LL v_n(x) & = \Delta v_n(x)-V(x)v_n(x) \\
& =  \frac{1}{\mu(x)}\sum_{y\in G}\omega(x,y)(v_n(y)-v_n(x))-V(x)v_n(x) = 0\quad
\forall\,\,x\in B_n(p),\,\,n\in\N,
\end{split}
\end{equation}
and since the sum which defines the Laplacian $\Delta$ is
actually a \emph{finite sum}
(recall that the graph $G$ is lo\-cal\-ly finite), by letting $n\to+\infty$ in \eqref{eq:LLvzero} we
readily obtain
$$\LL v(x) = 0\quad\forall\,\,x\in \bigcup_{n = 1}^{+\infty} B_n(p) =  G,$$
and thus $v$ is a solution of problem \eqref{problema}. Finally,
reminding that $0\leq v\leq 1$ on $G$ and using the Strong Maximum Principle
in Lemma \ref{lem:SMP}, we conclude that $v\in\mathfrak{F}$
is a solution of problem \eqref{problema} also fulfilling \eqref{e1g}.
Hence, we are left to prove the claimed (i)-(ii).
\vspace*{0.1cm}

Let us show (i). We first observe that, since $\LL v_n = 0$ in $B_n(p)$
 and since $v_n = (u-M)_+\geq 0$ in $G\setminus B_n(p)$, from the
 Weak Maximum Principle in Lemma \ref{lem:WMP} we infer that $v_n\geq 0$ in $G$.
 On the other hand, since the constant function $\zeta\equiv 1$ satisfies
 $$
 \begin{cases}
 \LL\zeta = -V(x)\leq 0 & \text{in $B_n(p)$}, \\
 \zeta \geq (u-M)_+ & \text{in $G\setminus B_n(p)$},
 \end{cases}$$
 by applying the Weak Maximum Principle in Lemma \ref{lem:WMP} to the function $\zeta-v_n$
 we immediately conclude that $v_n\leq \zeta\equiv 1$ pointwise on $G$.
 \vspace*{0.1cm}

Now, let us prove (ii). First of all we observe that, since $\zeta = u-M$ satisfies
 \[
\begin{cases}
\LL\zeta = VM\geq 0 & \text{in $B_n(p)$}, \\
\zeta = u-M\leq (u-M)_+ & \text{in $G\setminus B_n(p)$},
\end{cases}
\]
by applying the Weak Maximum Principle in Lemma \ref{lem:WMP} to the function $v_{n+1}-\zeta$ we derive that
$v_{n+1}\geq \zeta = u-M$ in $G$; thus, since we already know that $v_{n+1}\geq 0$, we obtain
\begin{equation} \label{eq:tousevnvnp1}
 v_{n+1} \geq (u-M)_+\quad\text{pointwise on $G$}.
\end{equation}
Owing to \eqref{eq:tousevnvnp1},
and applying once again the Weak Maximum Principle in Lemma \ref{lem:WMP} to the function
$v_{n+1}-v_{n}$ with $\Omega = B_n(p)$, we then conclude that $v_n\leq v_{n+1}$ in $G$.
\vspace*{0.1cm}

\noindent This ends the proof.
 \end{proof}

\section{A useful apriori estimate}\label{apriori}\setcounter{equation}{0}
Now we have established Proposition \ref{prop2}, we turn to prove
an \emph{apriori} estimate for non\-ne\-ga\-tive
and bounded solutions to \eqref{problema} which will play a key role in the proof of Theorem \ref{teo1}. Throughout what follows, we set
\begin{equation}\label{S}
S:=G\times [0, +\infty);
\end{equation}
and, for any given $T > 0$,
\[S_T:=G\times [0, T]\,.\]
\begin{proposition}\label{lemma1}
Let assumption \eqref{e7f} be in force. Let $u$ be a solution of equation \eqref{problema} such that $0\leq u\leq 1$. Moreover, let $T>0$ and define
\[v(x,t):=e^{t}u(x) - 1 \quad \text{for all }\; (x,t)\in S_T\,.\]
Finally, let $\eta\in \mathfrak F, $ and $\xi:S_T\to \mathbb R$ be such that
\begin{align}
&\mathrm{(i)}\quad \eta\ge0, \; \operatorname{supp} \eta \text{ is finite };\label{e9f}\\
&\mathrm{(ii)}\quad \xi(x, \cdot)\in C^1([0, T)),\,\, e^{\xi(x, \cdot)}\in C^1([0, T]) \quad \text{for any}\,\, x\in G;\label{e9bf}\\
&\mathrm{(iii)}\quad [\eta^2(y)-\eta^2(x)][e^{\xi(y,t)}-e^{\xi(x,t)}]\ge0 \quad \text{for all}\,\, x,y\in G,\,\, x\sim y,\,\, t\in [0, T].\label{e10f}
\end{align}
Then

\begin{equation}\label{eq31}
\begin{aligned}
&\sum_{x\in G} V(x) \eta^2(x) v^2_+(x, T)e^{\xi(x. T)}\mu(x) - \sum_{x\in G} V(x) \eta^2(x) v^2_+(x, 0)e^{\xi(x. 0)}\mu(x) \\ &
\quad \leq \int_0^T \sum_{x\in G} v_+^2(x,t)\eta^2(x)e^{\xi(x,t)}\left\{V(x) \xi_t(x,t)\mu(x) +\frac 12\sum_{y\in G}\omega(x,y)[1-e^{\xi(y,t)-\xi(x,t)}]^2\right\}dt\\
&\quad\quad +2\int_0^T\sum_{x,y\in G} v^2_+(x,t)e^{\xi(y,t)}[\eta(y)-\eta(x)]^2\omega(x,y)\,dt\,.
\end{aligned}
\end{equation}
\end{proposition}
To prove Proposition \ref{lemma1} we need an auxiliary result
(of independent interest) which shows that, if
$v$ is any solution of the parabolic equation
\begin{equation} \label{eq31g}
V\partial_t u-\Delta u=0\quad\text{in $G\times[0,T]$},
\end{equation}
then
$
v_+:=\max\{v;\,0\}
$
is a subsolution of the same equation (in a suitable sense).
\medskip

This is the content of the following
\begin{lemma}\label{lem4}
Let $v:S_T\to \R$ be a solution of equation \eqref{eq31g} such that the map $t\mapsto v(x,t)$ is $C^1([0,T])$ for any $x\in G$. Then, for any fixed $x\in G$,
$$
V\partial_ t v_+(x,t)-\Delta v_+(x,t)\le 0\quad \text{for a.e.}\,\,\,t\in [0,T].
$$
\end{lemma}

\begin{proof}
We separately consider three cases.

\smallskip

\noindent \textit{Case 1.} Let $(x,t)\in G\times(0,T)$ be such that
\begin{equation}\label{eq1111}
v(x,t)>0\,.
\end{equation}
Then, by the continuity of the map $t\mapsto v(x,t)$, there exists $\delta>0$ such that $v(x,\tau)>0$ for any $\tau\in(t-\delta, t+\delta)$.
Therefore, $v(x,\tau)=v_+(x,\tau)$ for any $\tau\in(t-\delta, t+\delta)$. Consequently, since the map $t\mapsto v(x,t)$ is $C^1([0,T])$, we also have that \begin{equation}\label{eq1111b}
\partial_t v(x,t)=\partial_t v_+(x,t).
\end{equation}
Now, due to the fact that $v_+(y,t)\ge v(y,t)$ for any $y\in G$, by \eqref{eq1111}, \eqref{eq1111b} and \eqref{eq31g}, we have
$$
\begin{aligned}
V\partial_t v_+(x,t)-\Delta v_+(x,t)&=V\partial_t v(x,t)-\frac{1}{\mu(x)}\sum_{y\in G}[ v_+(y,t)-v_+(x,t)]\, \omega(x,y)\\
&\le V\partial_t v(x,t)-\frac{1}{\mu(x)}\sum_{y\in G}[ v(y,t)-v(x,t)]\, \omega(x,y)=0\,.
\end{aligned}
$$

\smallskip

\noindent \textit{Case 2.} Let $(x,t)\in G\times(0,T)$ be such that
\begin{equation}\label{eq1112}
v(x,t)<0\,.
\end{equation}
Then, by the continuity of the map $t\mapsto v(x,t)$, there exists $\delta>0$ such that $v(x,\tau)<0$ for any $\tau\in(t-\delta, t+\delta)$.
Therefore, $v_+(x,\tau)=0$ for any $\tau\in(t-\delta, t+\delta)$. Consequently, since the map $t\mapsto v(x,t)$ is $C^1([0,T])$, we also have that
\begin{equation}\label{eq1112b}
\partial_t v_+(x,t)=0.
\end{equation}
Now, due to the fact that $v_+(y,t)\ge 0$ for any $y\in G$, by \eqref{eq1112} and \eqref{eq1112b}, we have
$$
\begin{aligned}
V\partial_t v_+(x,t)-\Delta v_+(x,t)&=-\frac{1}{\mu(x)}\sum_{y\in G}[ v_+(y,t)-v_+(x,t)] \,\omega(x,y)\\
&=-\frac{1}{\mu(x)}\sum_{y\in G}[ v_+(y,t)] \omega(x,y)\le 0\,.
\end{aligned}
$$

\smallskip

\noindent \textit{Case 3.} Let $x\in G$ fixed. We consider the set
$\mathcal{U}_x=\{t\in(0,T):v(x,t)=0\}$. We have that
\begin{equation}\label{eq1113}
\partial_t v(x,t)=0\quad \text{for a.e.}\,\,\,t\in\mathcal{U}_x.
\end{equation}
Then we can find $\mathcal{V}_x\subset\mathcal{U}_x$ such that
\begin{equation}\label{eq1113b}
\begin{aligned}
\text{$\lambda(\mathcal{V}_x)=\lambda(\mathcal{U}_x)$ and
$\partial_t v(x,t)=0$ for every $t\in \mathcal{V}_x$}
\end{aligned}
\end{equation}
(here, $\lambda$ denotes the one-dimensional Lebesgue measure). Thus, since $v_+(y,t)\ge 0$ for
pointwise on $G\times(0,T)$, by \eqref{eq1113} and \eqref{eq1113b}, for any $t\in\mathcal{V}_x$, we have
$$
\begin{aligned}
V\partial_t v_+(x,t)-\Delta v_+(x,t)&=-\frac{1}{\mu(x)}\sum_{y\in G}[ v_+(y,t)-v_+(x,t)]\, \omega(x,y)\\
&=-\frac{1}{\mu(x)}\sum_{y\in G}[ v_+(y,t)]\, \omega(x,y)\le 0\,.
\end{aligned}
$$
Hence the thesis follows from the combination of the three cases.
\end{proof}

\begin{proof}[Proof of Proposition \ref{lemma1}]

We first observe that, due to Lemma \ref{lem4}, for any $x\in G$,
\begin{equation}\label{eq32}
V(x)\,\partial_t v_+(x,t)-\Delta v_+(x,t)\le 0 \quad \text{for a.e.}\,\,t\in[0,T].
\end{equation}
We multiply \eqref{eq32} by
the function $(x,t)\mapsto v_+(x,t)\,\eta^2(x)e^{\xi(x,t)}\mu(x)$ and we integrate on
the time interval $[0,T]$. Thus we get
$$
\begin{aligned}
\int_0^T V(x)\, &\partial_t v_+(x,t) v_+(x,t) \,\eta^2(x)e^{\xi(x, t)}\mu(x)\,dt\\
&\le \int_0^T \Delta(v_+(x,t))\,v_+(x,t)\,\eta^2(x)e^{\xi(x,t)}\mu(x)\,dt.
\end{aligned}
$$
We sum over $x\in G$, so
\begin{equation*}
\begin{split}
\begin{aligned}
\sum_{x\in G}\int_0^T V(x)\,&\partial_t v_+(x,t)  v_+(x,t) \,\eta^2(x)e^{\xi(x,t)}\mu(x)\,dt\\
&\le\sum_{x\in G} \int_0^T \Delta(v_+(x,t))\,v_+(x,t)\,\eta^2(x)e^{\xi(x,t)}\mu(x)\,dt.
\end{aligned}
\end{split}
\end{equation*}
Notice that, since $\eta(x)$ is finitely supported, the series in the latter inequality are actually finite sums. Therefore,
we obtain
\begin{equation}\label{eq33}
\begin{aligned}
\int_0^T\sum_{x\in G}\,V(x)\,&\partial_t v_+(x,t)  v_+(x,t) \,\eta^2(x)e^{\xi(x,t)}\mu(x)\,dt\\
&\le \int_0^T\sum_{x\in G}\Delta(v_+(x,t))\,v_+(x,t)\,\eta^2(x)e^{\xi(x,t)}\mu(x)\,dt.
\end{aligned}
\end{equation}
Set
\begin{equation*}
I:=\sum_{x\in G}\Delta(v_+(x,t))\,v_+(x,t)\,\eta^2(x)e^{\xi(x,t)}\mu(x)\,.
\end{equation*}
Then, by \eqref{e2f} and \eqref{e4f},
\begin{equation}\label{eq34}
\begin{aligned}
I&=-\frac 12 \sum_{x,y\in G}\left(\nabla_{xy} v_+\right) \nabla_{xy}\left[v_+\eta^2e^{\xi}\right]\omega(x,y)\\
&=-\frac 12 \sum_{x,y\in G}\left[\nabla_{xy}v_+\right]^2\eta^2(y)e^{\xi(y,t)}\,\omega(x,y)\\
&\quad-\frac 12 \sum_{x,y\in G}v_+(x,t)\eta^2(x)\left(\nabla_{xy}v_+\right)\left(\nabla_{xy}e^{\xi}\right)\,\omega(x,y)\\
&\quad-\frac 12 \sum_{x,y\in G}v_+(x,t)e^{\xi(y,t)}\left(\nabla_{xy}v_+\right)\left(\nabla_{xy}\eta^2\right)\,\omega(x,y)=:J_1+J_2+J_3\,.
\end{aligned}
\end{equation}
In view of \eqref{e2f}, we obviously have
\begin{equation}\label{e11f}
\dif \eta^2=[\eta(y)+\eta(x)][\eta(y)-\eta(x)] \quad \text{ for all }\, x,y\in G\,.
\end{equation}
By Young's inequality with exponent $2$, for any $\delta_1>0$, we have, for every $t\in [0,T]$,
\begin{equation}\label{eq38}
\begin{aligned}
J_2&=-\frac 12 \sum_{x,y\in G}v_+(x,t)\eta^2(x)(\nabla_{xy}v_+)(\nabla_{xy}e^{\xi})\,\omega(x,y)\\
&=-\frac 12 \sum_{x,y\in G}v_+(y,t)\eta^2(y)(\nabla_{xy}v_+)\left[e^{\xi(y,t)}-e^{\xi(x,t)}\right]\,\omega(x,y)\\
&=-\frac 12 \sum_{x,y\in G}v_+(y,t)\eta^2(y)(\nabla_{xy}v_+) e^{\xi(y,t)}\left[1-e^{\xi(x,t)-\xi(y,t)}\right]\,\omega(x,y)\\
&\le \frac{\delta_1}{4}\sum_{x,y\in G}\left[\nabla_{xy}v_+\right]^2\eta^2(y)e^{\xi(y,t)}\,\omega(x,y)\\
&\quad + \frac1{4\delta_1}\sum_{x,y\in G}v_+^2(y,t)\left[1-e^{\xi(x,t)-\xi(y,t)}\right]^2\eta^2(y)e^{\xi(y,t)}\,\omega(x,y)\\
&= \frac{\delta_1}{4}\sum_{x,y\in G}\left[\nabla_{xy}v_+\right]^2\eta^2(y)e^{\xi(y,t)}\,\omega(x,y)\\
&\quad + \frac1{4\delta_1}\sum_{x,y\in G}v_+^2(x,t)\left[1-e^{\xi(y,t)-\xi(x,t)}\right]^2\eta^2(x)e^{\xi(x,t)}\,\omega(x,y).\\
\end{aligned}
\end{equation}
Similarly, due to \eqref{e11f}, by Young's inequality with exponent $2$, we have, for every $\delta_2>0$ and for every $t\in[0,T]$,
\begin{equation}\label{eq35}
\begin{aligned}
J_3&= -\frac 12 \sum_{x,y\in G}v_+(x,t)e^{\xi(y,t)}(\nabla_{xy}v_+)\left[\eta(y)+\eta(x)\right]\left[\eta(y)-\eta(x)\right]\,\omega(x,y)\\
&\le \frac{\delta_2}{4} \sum_{x,y\in G}e^{\xi(y,t)}\left[\nabla_{xy}v_+\right]^2\left[\eta(y)+\eta(x)\right]^2\,\omega(x,y)\\
&\,\,\,\,\,+\frac 1{4\delta_2} \sum_{x,y\in G}e^{\xi(y,t)} v_+^2(x,t)\left[\eta(y)-\eta(x)\right]^2\,\omega(x,y)\\
&\le \frac{\delta_2}{2} \sum_{x,y\in G}e^{\xi(y,t)}\left[\nabla_{xy}v_+\right]^2 \left[\eta^2(y)+\eta^2(x)\right]\,\omega(x,y)\\
&\,\,\,\,\,+\frac 1{4\delta_2} \sum_{x,y\in G}e^{\xi(y,t)} v_+^2(x,t)\left[\eta(y)-\eta(x)\right]^2\,\omega(x,y)\\
&= \frac{\delta_2}{4} \sum_{x,y\in G}\left[e^{\xi(y,t)}+e^{\xi(x,t)}\right]\left[\nabla_{xy}v_+\right]^2\left[\eta^2(y)+\eta^2(x)\right]\,\omega(x,y)\\
&\,\,\,\,\,+\frac 1{4\delta_2} \sum_{x,y\in G}e^{\xi(y,t)} v_+^2(x,t)\left[\eta(y)-\eta(x)\right]^2\,\omega(x,y)\,.\\
\end{aligned}
\end{equation}
From \eqref{e10f} we can easily infer that
\begin{equation}\label{eq36}
\left[\eta^2(y)+\eta^2(x)\right]\left[e^{\xi(y,t)}+e^{\xi(x,t)}\right]\le 2\left\{\eta^2(y)e^{\xi(y,t)}+\eta^2(x)e^{\xi(x,t)}\right\}\,.
\end{equation}
By using \eqref{eq35} and \eqref{eq36},
\begin{equation}\label{eq37}
\begin{aligned}
J_3&\le \frac{\delta_2}{2} \sum_{x,y\in G}\left[\nabla_{xy}v_+\right]^2\left\{\eta^2(y)e^{\xi(y,t)}+\eta^2(x)e^{\xi(x,t)}\right\}\,\omega(x,y)\\
&\,\,\,\,\,\,\,+\frac 1{4\delta_2} \sum_{x,y\in G}e^{\xi(y,t)} v_+^2(x,t)\left[\eta(y)-\eta(x)\right]^2\,\omega(x,y)\\
&=\delta_2\sum_{x,y\in G}\left[\nabla_{xy}v_+\right]^2\eta^2(y)e^{\xi(y,t)}\,\omega(x,y)\\
&\,\,\,\,\,\,\,+\frac 1{4\delta_2} \sum_{x,y\in G}e^{\xi(y,t)} v_+^2(x,t)\left[\eta(y)-\eta(x)\right]^2\,\omega(x,y)\,.\\
\end{aligned}
\end{equation}
By combining \eqref{eq34}, \eqref{eq38} and \eqref{eq37} we get, for every $t\in[0,T]$,
\begin{equation}\label{eq39}
\begin{aligned}
\sum_{x\in G}&\Delta v_+(x,t)v_+(x,t),\eta^2(x)e^{\xi(x,t)}\mu(x)\\
&\le -\frac 12 \sum_{x,y\in G}[\nabla_{xy}v_+]^2\eta^2(y)e^{\xi(y,t)}\omega(x,y) \\
&\quad + \frac{\delta_1}{4}\sum_{x,y\in G}[\nabla_{xy}v_+]^2\eta^2(y)e^{\xi(y,t)}\omega(x,y) \\
&\quad+ \frac1{4\delta_1}\sum_{x,y\in G}v_+^2(x,t)\left[1-e^{\xi(y,t)-\xi(x,t)}\right]^2\eta^2(x)e^{\xi(x,t)}\,\omega(x,y) \\
&\quad + \delta_2\sum_{x,y\in G}\left[\nabla_{xy}v_+\right]^2\eta^2(y)e^{\xi(y,t)}\,\omega(x,y) \\
&\quad+\frac 1{4\delta_2} \sum_{x,y\in G} v_+^2(x,t)\left[\eta(y)-\eta(x)\right]^2e^{\xi(y,t)}\,\omega(x,y)\,.
\end{aligned}
\end{equation}
We now choose $\delta_1=1$ and $\delta_2=\frac 14$ in such a way that $-\frac 12+ \frac{\delta_1}{4}+ \delta_2=0$. Consequently, \eqref{eq39} yields, for \textcolor{red}{every} $t\in[0,T]$,
\begin{equation}\label{eq310}
\begin{aligned}
I&=\sum_{x\in G}\Delta v_+(x,t)v_+(x,t)\,\eta^2(x)e^{\xi(x,t)}\mu(x)\\
&\le \frac1{4}\sum_{x,y\in G}v_+^2(x,t)\left[1-e^{\xi(y,t)-\xi(x,t)}\right]^2\eta^2(x)e^{\xi(x,t)}\,\omega(x,y) \\
&\quad + \sum_{x,y\in G}v_+^2(x,t)\left[\eta(y)-\eta(x)\right]^2e^{\xi(y,t)}\,\omega(x,y)\,.
\end{aligned}
\end{equation}
By substituting \eqref{eq310} into \eqref{eq33}, due to the linearity of the integral operation, we get
\begin{equation}\label{eq311}
\begin{aligned}
\int_0^T \sum_{x\in G}\,&V(x)\,\partial_t v_+(x,t)\,v_+(x,t)\eta^2(x)e^{\xi(x,t)}\mu(x)\,dt\\
&\le \int_0^T\frac1{4}\sum_{x,y\in G}v_+^2(x,t)\left[1-e^{\xi(y,t)-\xi(x,t)}\right]^2\eta^2(x)e^{\xi(x,t)}\,\omega(x,y)\,dt \\
&\quad + \int_0^T\sum_{x,y\in G}v_+^2(x,t)[\eta(y)-\eta(x)]^2e^{\xi(y,t)}\omega(x,y)\,dt.
\end{aligned}
\end{equation}
We now consider the right hand side of \eqref{eq311}. Thus we have
\begin{equation}\label{eq314b}
\begin{aligned}
\int_0^T \sum_{x\in G}\,&V(x)\,\partial_t v_+(x,t)\,v_+(x,t)\eta^2(x)e^{\xi(x,t)}\mu(x)\,dt\\
&=\sum_{x\in G}\,V(x)\eta^2(x)\mu(x)\int_0^T(v_+)_t(x,t)\,v_+(x,t)e^{\xi(x,t)}\,dt\\
&=\frac 12\sum_{x\in G}\,V(x)\eta^2(x)\mu(x)\int_0^T(v_+^2)_t(x,t)\,e^{\xi(x,t)}\,dt\\
&=\frac 12\sum_{x\in G}\,V(x)\eta^2(x)v_+^2(x,T)\,e^{\xi(x,T)}\,\mu(x)\\
&\quad -\frac 12\sum_{x\in G}\,V(x)\eta^2(x)v_+^2(x,0)\,e^{\xi(x,0)}\,\mu(x)\\
&\quad -\frac 12 \int_0^T\sum_{x\in G}\,V(x)\eta^2(x)v_+^2(x,t)\xi_t(x,t)\,e^{\xi(x,t)}\,\mu(x)\,dt\,.
\end{aligned}
\end{equation}
By combining together \eqref{eq311} and \eqref{eq314b}, we then obtain
\begin{equation*}
\begin{aligned}
& \sum_{x\in G}\,V(x)\eta^2(x)v_+^2(x,T)\,e^{\xi(x,T)}\,\mu(x) -\sum_{x\in G}\,V(x)\eta^2(x)v_+^2(x,0)\,e^{\xi(x,0)}\,\mu(x)\\
&\quad \le\int_0^T\sum_{x\in G}\,\eta^2(x)v_+^2(x,t)\,e^{\xi(x,t)}\left\{V(x)\xi_t(x,t)\mu(x)+\frac 12 \sum_{y\in G}\left[1-e^{\xi(y,t)-\xi(x,t)}\right]^2\,\omega(x,y)\right\}\,dt\\
&\quad\quad + 2\int_0^T\sum_{x,y\in G}v_+^2(x,t)[\eta(y)-\eta(x)]^2e^{\xi(y,t)}\omega(x,y)\,dt.
\end{aligned}
\end{equation*}
This is precisely the \eqref{eq31}, and the proof is complete.
\end{proof}

\section{Some distinguished test functions}\label{test}\setcounter{equation}{0}
Let us now prove the existence of suitable test functions $\xi$ and $\eta$ which
are
\emph{admissible} in \eqref{eq31} and which satisfy
some \emph{ad-hoc} properties.
\vspace*{0.1cm}

Throughout the sequel, we let $x_0\in G$ and $R_0 > 0$
be as in assumption \eqref{e12f}; furthermore, we choose parameters
$M>0$, $T>0$, $\beta\in (0, 1]$, $\lambda>1$ and
\begin{equation}\label{e3t}
r\geq 2s+R_0,
\end{equation}
where $s$ is defined in \eqref{e14f}.
Setting
$$\text{$\mathbf{d}(x) := d(x,x_0)$ for every $x\in\ G$},$$
we then define
\begin{equation}\label{eq313}
\begin{gathered}
\xi(x,t):=- M \frac{\rho(x)}{\lambda T-t} \quad \text{for any}\,\, x\in G,\,\, t\in [0, T), \\
\text{where $\rho(x):= \max\{\mathbf{d}^\beta(x),r^\beta\} =
\begin{cases}
\,\,\,\,r^{\beta}  &\text{if} \,\,\, \mathbf{d}(x)\leq r,\\
\mathbf{d}^{\beta}(x) & \text{if} \,\,\, \mathbf{d}(x)>r\,,
\end{cases}$}
\end{gathered}
\end{equation}
Concerning the function $\xi$, we have the following key lemma.
\begin{lemma}\label{lemma2}
Let hypotheses \eqref{e1f} and \eqref{e7f} be fulfilled; suppose that \eqref{e12f} is satisfied with $\alpha\in [0, 1]$. Let $\xi$ be the function defined in \eqref{eq313}.
Then
\begin{equation}\label{eq314}
V(x)\,\xi_t(x,t)\mu(x)+\frac 12 \sum_{y\in G} \omega(x,y)\left[1-e^{\xi(y,t)-\xi(x,t)}\right]^2
\leq 0 \quad \text{for all}\,\,\, x\in G,\,\, t\in (0, T),
\end{equation}
provided that $M=T>0$ is sufficiently small and $\beta=\alpha$.
\end{lemma}
\begin{proof}
 To ease the readability, we split the proof into two steps.
 \medskip

 \noindent \textsc{Step I:} In this first step we prove the following estimate
 \begin{equation} \label{eq:estimrhoStepI}
  |\rho(x)-\rho(y)|\leq \beta R_0^{\beta-1}d(x,y) \qquad
  \text{for every $x,y\in G$ with $x\sim y$}.
 \end{equation}
 To this end, it is useful to distinguish two cases.
 \begin{itemize}
  \item[(i)] $x\in B_{r-s}(x_0)$. In this case, by triangle's inequality and \eqref{e14f} we have
 $$\text{$\mathbf{d}(y)< r$ for every $y\in G,\,y\sim x$};$$
 as a consequence, from the very definition of $\rho$ we derive
 $$|\rho(x)-\rho(y)| = r^\beta - r^\beta = 0,$$
 and this trivially implies \eqref{eq:estimrhoStepI}.
 \medskip

 \item[(ii)] $x\in G\setminus B_{r-s}(x_0)$.
 In this case we first notice that, since the function
 $$\R\ni t\mapsto\max\{t,r^\beta\}$$
 is Lipschitz-continuous with Lipschitz constant $L = 1$,
 by the Mean Value Theorem and again the triangle inequality we can write
 \begin{equation} \label{eq:estimrhoCaseiigeneral}
  \begin{split}
   |\rho(x)-\rho(y)| & \leq |\mathbf{d}^\beta(x)-\mathbf{d}^\beta(y)|
   \leq \beta\sigma^{\beta-1}|\mathbf{d}^\beta(x)-\mathbf{d}^\beta(y)| \\
   & \leq \beta\sigma^{\beta-1}d(x,y)\qquad\forall\,\,y\in G,
  \end{split}
 \end{equation}
 where $\sigma \geq 0$ is a suitable point between $\mathbf{d}(x)$ and $\mathbf{d}(y)$.
 On the other hand, since we as\-su\-ming that $x\notin B_{r-s}(x_0)$
 (hence, $\mathbf{d}(x)\geq r-s$),
 by \eqref{e14f} we have
 \begin{equation} \label{eq:dygeq}
  \mathbf{d}(y)\geq r-2s\geq R_0\quad\forall\,\,y\in G,\,y\sim x.
  \end{equation}
 Recalling that $\beta \leq 1$, from \eqref{eq:estimrhoCaseiigeneral}-\eqref{eq:dygeq}
 we immediately obtain
 $$|\rho(x)-\rho(y)|\leq \beta R_0^{\beta-1}d(x,y)\quad\forall\,\,y\in G,\,y\sim x,$$
 which is exactly the desired \eqref{eq:estimrhoStepI}.
 \end{itemize}
 \medskip

 \noindent\textsc{Step II:} In this second step we establish \eqref{eq314}.
 To begin with, we point out that
 $$(e^a-1)^2\leq a^2e^{2|a|}\quad\forall\,\,a\in\R;$$
 this inequality, together with \eqref{eq:estimrhoStepI}, allows us to write
 \begin{align*}
  \big(1-e^{\xi(y,t)-\xi(x,t)}\big)^2 &
  \leq |\xi(y,t)-\xi(x,t)|^2 e^{2|\xi(y,t)-\xi(x,t)|} \\
  & = \frac{M^2}{(\lambda T-t)^2}|\rho(y)-\rho(x)|^2
  e^{\frac{2M}{\lambda T-t}|\rho(y)-\rho(x)|} \\
  & \leq
  \frac{\beta^2R_0^{2\beta-2}M^2}{(\lambda T-t)^2}
  e^{\frac{2\beta R_0^{\beta -1} M d(x,y)}{T(\lambda-1)}}d(x,y)^2,
 \end{align*}
 and this estimate holds for every $x,y\in G$ and every $t\in (0,T)$.
 Then, using \eqref{e14f}
 and recalling that $d$ is \emph{intrinsic}
 (hence, \eqref{e1f} holds),
 for every $x\in G$ and $t\in (0,T)$ we obtain
 \begin{equation}  \label{eq:estimsupersolI}
 \begin{split}
   & V(x)\,\xi_t(x,t)\mu(x)+\frac 12 \sum_{y\in G} \omega(x,y)\left[1-e^{\xi(y,t)-\xi(x,t)}\right]^2
   \\
   & \qquad \leq -\frac{M\rho(x)}{(\lambda T-t)^2}
   V(x)\mu(x)+\frac{\beta^2R_0^{2\beta-2}M^2}{2(\lambda T-t)^2}
   \sum_{y\in G}e^{\frac{2\beta R_0^{\beta -1} M d(x,y)}{T(\lambda-1)}}\omega(x,y) d(x,y)^2 \\
   & \qquad \leq
   -\frac{M\rho(x)}{(\lambda T-t)^2}
   V(x)\mu(x)+\frac{\beta^2R_0^{2\beta-2}M^2}{2(\lambda T-t)^2}
   e^{\frac{2\beta R_0^{\beta -1} M s}{T(\lambda-1)}}
   \sum_{y\in G}\omega(x,y) d(x,y)^2 \\
   & \qquad = \frac{M\mu(x)}{(\lambda T-t)^2}\left\{-V(x)\rho(x)+
    \frac{M}{2}\beta^2R_0^{2\beta-2}e^{\frac{2\beta R_0^{\beta -1} M s}{T(\lambda-1)}}\right\}.
  \end{split}
 \end{equation}
 To proceed further, we now \emph{fix} $\beta = \alpha$ and we
 exploit assumption \eqref{e12f}:
 taking into account the piecewise definition of $\rho$,
 see \eqref{eq313}, it
 is easy to recognize that
 \begin{equation} \label{eq:assVgeq}
  V(x)\rho(x)\geq c_0\,\mathbf{d}(x)^{-\alpha}\rho(x)
 = c_0\,\mathbf{d}(x)^{-\alpha}\cdot\max\{\mathbf{d}^\alpha(x),r^\alpha\}\geq c_0\quad
 \forall\,\,x\in G;
 \end{equation}
 as a consequence, by combining \eqref{eq:estimsupersolI}-\eqref{eq:assVgeq} we conclude that
 \begin{align*}
 & V(x)\,\xi_t(x,t)\mu(x)+\frac 12 \sum_{y\in G} \omega(x,y)\left[1-e^{\xi(y,t)-\xi(x,t)}\right]^2 \\
 &\qquad \leq
 \frac{M\mu(x)}{(\lambda T-t)^2}\left\{-c_0+
    \frac{M}{2}\beta^2R_0^{2\beta-2}e^{\frac{2\beta R_0^{\beta -1} M s}{T(\lambda-1)}}\right\}\leq 0,
 \end{align*}
 provided that
\begin{equation*}
M=T\quad\text{and}\quad
 0<M\leq  \frac{2 c_0 R_0^{2-2\beta} e^{-\frac{2}{\lambda-1}\beta R_0^{\beta-1}s}}{ \beta^2}\,.
\end{equation*}
This ends the proof.
\end{proof}
Now that we have proved Lemma \ref{lemma2}, we turn to prove the existence
of a suitable `cut-off' function $\eta$. To this end, taking for fixed all the notation
introduced so far, we choose
$$r_1 > 2r+8s$$
and we define the function
\begin{equation}\label{eq316}
\eta(x):=\min\left\{\frac{2\left[r_1-s-\mathbf{d}(x)\right]_+}{r_1},\,\,1\right\}\, \quad \text{ for any }\; x\in G\,.
\end{equation}
Owing to \cite[Lemma 5.2]{MP} (with the choice $\delta = 1/2$), we obtain the following result.
\begin{lemma}\label{lemma3}
Let assumptions \eqref{e1f}-\eqref{e7f} be satisfied. Then, the function $\eta$ defined in \eqref{eq316}
satisfies the following properties:
\begin{align}
& \mathrm{(i)}\quad |\dif \eta| \leq \frac 2{ r_1} d(x,y)\chi_{\big\{\frac{r_1}{2}-2s
 \le \mathbf{d}(x)\le r_1\big\}} \quad\quad \text{ for any }\; x\in G; \nonumber \\
& \mathrm{(ii)}\quad \sum_{y\in G} \left(\dif \eta\right)^2 \omega(x,y)\le \frac{4}{ (r_1)^2}\mu(x)
\chi_{\big\{\frac{r_1}{2}-2s
 \le \mathbf{d}(x)\le r_1\big\}} \quad \text{ for any }\; x\in G. \label{eq317}
\end{align}
\end{lemma}

\section{Proof of Theorem \ref{teo1}}\label{proofthm}\setcounter{equation}{0}
Due to the results established in Sections \ref{auxiliary}, \ref{apriori} and \ref{test}, we are ready to provide the proof of Theorem \ref{teo1}. In what follows,
we take for fixed all the notation introduced so far.

\begin{proof}[Proof of Theorem \ref{teo1}]
By contradiction, suppose that there exists a non-trivial bounded solution of equation \eqref{problema}. Then, due to Proposition \ref{prop2}, we know that there exists a solution $u$ to to same equation \eqref{problema} such that
  \begin{equation} \label{eq:uzeropositive}
   0 < u \leq 1\quad\text{pointwise in $G$}.
  \end{equation}
Let us now define $v(x,t) := e^{t}u(x)-1$, for any $(x,t)\in S$, for $S$ as in \eqref{S}.
We want to show that
 \begin{equation} \label{eq:claimcentrale}
   v(x,t)\leq 0\quad\text{for every $x\in\mathrm{supp}(V)$ and $t > 0$}.
  \end{equation}

  \noindent To do so, let us fix $T > 0$ (to be chosen conveniently small in a moment), and we arbitrarily choose $r >  2s+R_0$, for $R_0 > 0$ and $0<s<+\infty$ as in \eqref{e12f} and \eqref{e14f}, respectively. Let $\xi$ be as in \eqref{eq313} with $\beta=\alpha$, $\lambda>1$ and with $M=T$ chosen as in Lemma \ref{lemma2}.
  Moreover, let us fix $r_1 > 0$ in such a way that
  \begin{equation} \label{eq:choiceR1R}
   r_1\geq 2r+8s,
  \end{equation}
  and let $\eta$ be as in Lemma \ref{lemma3}.
 Now, we observe that $\eta$ and $\xi$ obviously satisfy conditions \eqref{e9f} and \eqref{e9bf};  furthermore, also \eqref{e10f} is fulfilled, since both $\eta$ and $\xi(\cdot,t)$ are \emph{non-increasing functions} of $\mathbf{d}$. Therefore, form \eqref{eq31} we obtain
  \begin{equation} \label{eq:estimdastimarePartI}
   \begin{aligned}
   &\sum_{x\in G} V(x) \eta^2(x) v^2_+(x, T)e^{\xi(x. T)}\mu(x) - \sum_{x\in G} V(x) \eta^2(x) v^2_+(x, 0)e^{\xi(x. 0)}\mu(x) \\
&\quad
\leq \int_0^T \sum_{x\in G} v_+^2(x,t)\eta^2(x)e^{\xi(x,t)}\left\{V(x) \xi_t(x,t)\mu(x) -\frac 12\sum_{y\in G}\omega(x,y)[1-e^{\xi(y,t)-\xi(x,t)}]^2\right\}\,dt\\
&\qquad\qquad +2\int_0^T\sum_{x,y\in G} v^2_+(x,t)e^{\xi(y,t)}[\eta(y)-\eta(x)]^2\omega(x,y)\,dt.
\end{aligned}
\end{equation}
 On the other hand, by Lemma \ref{lemma2}
 there exists $T_0 = T_0(\lambda)> 0$, such that
 \begin{equation} \label{eq:disugLemma2zeta}
  V(x) \xi_t(x,t)\mu(x) -\frac 12\sum_{y\in G}\omega(x,y)[1-e^{\xi(y,t)-\xi(x,t)}]^2
   \leq 0\qquad\forall\,\,x\in G,\,t\in (0,T),
 \end{equation}
 provided that $T \leq T_0$. As a consequence, by combining \eqref{eq:estimdastimarePartI} and \eqref{eq:disugLemma2zeta}, we obtain
 \begin{equation} \label{eq:estimdastimare}
 \begin{split}
  & \sum_{x\in G} V(x) \eta^2(x) v^2_+(x, T)e^{\xi(x. T)}\mu(x) - \sum_{x\in G} V(x) \eta^2(x) v^2_+(x, 0)e^{\xi(x. 0)}\mu(x) \\
& \qquad \leq 2\int_0^T\sum_{x,y\in G} v^2_+(x,t)e^{\xi(y,t)}[\eta(y)-\eta(x)]^2\omega(x,y)\,dt.
\end{split}
 \end{equation}
 We then proceed by estimating both sides of \eqref{eq:estimdastimare}.
 \medskip

 \noindent-\,\,\emph{Estimate of the left-hand side.}
 First of all we observe that, owing to the definition of $\eta$ in
 \eqref{eq316}, we have $\eta\geq 0$ pointwise on $G$ and
 $\eta\equiv 1$ on $B_r(x_0) = \{x:\,\mathbf{d}(x) < r\}$.
Now we observe that, due to \eqref{eq:uzeropositive}, $v_+(x,0) = 0$ in $G$, moreover, $\rho(x)=\mathbf{d}(x)^{\alpha}$ in $B_r(x_0)$, hence we obtain
 \begin{equation} \label{eq:estimLHS}
   \begin{split}
    \sum_{x\in G} V(x) &\eta^2(x) v^2_+(x, T)e^{\xi(x. T)}\mu(x) - \sum_{x\in G} V(x) \eta^2(x) v^2_+(x, 0)e^{\xi(x. 0)}\mu(x)\\
    & = \sum_{x\in G}V(x)\eta(x)^2v^2_+(x,T) e^{\zeta(x,T)}\mu(x) \\
    & \geq \!\!\!\sum_{x\in B_r(x_0)}\!\!\!V(x)v^2_+(x,T) e^{\zeta(x,T)}\mu(x) \\
    &=\gamma_{r,\Lambda}\!\!\sum_{x\in B_r(x_0)}\!\!\!V(x)\,v^2_+(x,T)\mu(x),
   \end{split}
   \end{equation}
   where we have used the shorthand notation
   $\gamma_{r,\Lambda} = e^{-\frac{r^\alpha}{\lambda-1}} > 0$.
   \medskip

  \noindent -\,\,\emph{Estimate of the right-hand side.}
  We first observe that,
  using the definition of $\xi$ given in \eqref{eq313}, since $0\leq v_+(x,t)\leq e^t$ on $S$, for $S$ as in \eqref{S}, and by exploiting \eqref{eq317}, we have
\begin{align*}
    \int_0^T\sum_{x,y\in G}& v^2_+(x,t)e^{\xi(y,t)}[\eta(y)-\eta(x)]^2\omega(x,y)\,dt \\
    & = \int_0^T\sum_{x,y\in G} v^2_+(x,t)e^{\xi(y,t)}(\nabla_{xy}\eta)^2\omega(x,y)\,dt \\
    & = \int_0^T\sum_{x,y\in G} v^2_+(y,t)e^{\xi(x,t)}(\nabla_{xy}\eta)^2\omega(x,y)\,dt \\
    & \leq e^{2T}\int_0^T\sum_{x,y\in G}\exp\Big(-\frac{T\rho(x)}{\lambda T-t}\Big)
    (\nabla_{xy}\eta)^2\omega(x,y)\,dt \\
    & \leq Te^{2T}\sum_{x,y\in G}\exp\Big(-\frac{\rho(x)}{\lambda}\Big)
    (\nabla_{xy}\eta)^2\omega(x,y)\\
    & = \frac{4Te^{2T}}{(r_1)^2} \sum_{x\in G}\exp\Big(-\frac{\rho(x)}{\lambda}\Big)\chi_{\big\{\frac{r_1}{2}-2s\le d(x,x_0)\le r_1\big\}}(x)\mu(x) \\
    & \leq  \frac{4Te^{2T}}{(r_1)^2}\sum_{x\in G}\exp\Big(-\frac{\rho(x)}{\lambda}\Big)\chi_{\{r < d(x,x_0)\leq r_1\}}(x)\mu(x).
\end{align*}
Now, since $\rho(x) = \mathbf{d}^\alpha(x)$ when $\mathbf{d}(x) > r$, by using \eqref{eq:choiceR1R}, we obtain that
\begin{equation}\label{eq:estimRHS}
\begin{aligned}
    \int_0^T\sum_{x,y\in G} & v^2_+(x,t)e^{\xi(y,t)}[\eta(y)-\eta(x)]^2\omega(x,y)\,dt \\
    & \le \frac{4Te^{2T}}{(r_1)^2}\sum_{x\in G}\exp\Big(-\frac{\mathbf{d}^\alpha(x)}{\lambda}\Big)\chi_{\{r< d(x,x_0)\le r_1\}}(x)\mu(x) \\
    & = \frac{4Te^{2T}}{(r_1)^2}\sum_{x\in G}\exp\Big(-\frac{\mathbf{d}^\alpha(x)}{\lambda}\Big)\chi_{\{r< d(x,x_0)\le r_1\}}(x)\mu(x) \\
    & = \frac{c(T)}{(r_1)^2}\sum_{x\in G}\exp\Big(-\frac{\mathbf{d}^\alpha(x)}{\lambda}\Big)\chi_{\{r< d(x,x_0)\le r_1\}}(x)\mu(x).
  \end{aligned}
  \end{equation}
Then, by combining \eqref{eq:estimLHS} and \eqref{eq:estimRHS} with \eqref{eq:estimdastimare}, we obtain
   \begin{equation}  \label{eq:estimToConcludeFin}
   \begin{split}
    & \sum_{x\in B_r(x_0)}\!\!\!V(x)\,v^2_+(x,T)\mu(x)
    \leq \frac{c(T,\lambda,r)}{(r_1)^2}\sum_{x\in G}
    \exp\Big(-\frac{\mathbf{d}^\alpha(x)}{\lambda}\Big)\chi_{\{r< \mathbf{d}(x)\le r_1\}}(x)\mu(x),
  \end{split}
  \end{equation}
  where $c(T,\lambda,r) > 0$ is a constant only depending $T$, $\lambda$ and $r$.
Now, if $0<\Lambda <1$ is as in \eqref{eq:hpStrada3} and if we set
  $$
  \lambda = \frac{1}{\Lambda} > 1,
  $$
  then estimate \eqref{eq:estimToConcludeFin}  boils down to
  \begin{equation} \label{eq:topasslimite1}
   \begin{split}
   \sum_{x\in B_r(x_0)}\!\!\!V(x)\,v^2_+(x,T)\mu(x)
    & \leq \frac{c(T,\lambda,r)}{(r_1)^2}\sum_{x\in G}
    e^{-\Lambda \mathbf{d}^\alpha(x)}\chi_{\{r< \mathbf{d}(x)\le r_1\}}(x)\mu(x) \\
    & =  \frac{c(T,\lambda,r)}{(r_1)^2}\sum_{x\in G}
    e^{-\Lambda d^\alpha(x,x_0)}\chi_{\{r< \mathbf{d}(x)\le r_1\}}(x)\mu(x).
  \end{split}
  \end{equation}
  Furthermore, we recall that $r_1$ was arbitrarily fixed, hence by taking the limit as $r_1\to+\infty$ in \eqref{eq:topasslimite1} and by assumption \eqref{eq:hpStrada3}, we get
  \begin{equation*}
   \sum_{x\in B_r(x_0)}\!\!\!V(x)\,v^2_+(x,T)\mu(x) = 0.
  \end{equation*}
  Thus, we readily derive that
  \begin{equation*}
   v(x,T)\leq 0\quad\forall\,\,x\in B_r(x_0)\cap\mathrm{supp}(V),\,\,\,0\leq T\leq T_0(\lambda),
   \end{equation*}
  where we recall that the number $T_0$ depends on $\lambda$, which is by now fixed.
  From this, recalling also that $r \geq 2s+R_0$ was arbitrarily fixed, we
  then obtain
  \begin{equation}\label{e300f}
  v(x,T)\leq 0\quad\forall\,\,x\in\mathrm{supp}(V),\,\,\,0\leq T\leq T_0(\lambda).
  \end{equation}
  Now, let us introduce the {\it `shifted' function}
  $$v_1(x,t) = v(x,t+T_0(\lambda)).$$
  Clearly, it is still a solution of the parabolic problem \eqref{eq31g}. In addition,
  by
  \eqref{e300f}, $v_1(x,0)\leq 0$ for every  $x\in \mathrm{supp}(V)$.
By applying
  the very same argument exploited so far, we can infer that
  $$v_1(x,T)\leq 0\quad\forall\,\,x\in\mathrm{supp}(V),\,0\leq T\leq T_0(\lambda),$$
  hence
  \[v(x,T)\leq 0\quad\forall\,\,x\in\mathrm{supp}(V),\,\,\,0\leq T\leq 2 T_0(\lambda).\]
  By iterating this argument, and by using in a crucial way
  the fact that $T_0(\lambda) > 0$
  is a \emph{universal number remaining unchanged at any iteration}
  (as $\lambda > 0$ is fixed), we conclude that
  $$v(x,T)\leq 0\quad\forall\,\,x\in\mathrm{supp}(V),\,T \geq 0.$$
  This yelds \eqref{eq:claimcentrale}.
  \medskip

  We can now easily conclude the proof of the theorem, in fact, due to \eqref{eq:claimcentrale} and exploiting the definition of $v$, we have
  $$
  0 < u\leq e^{-t}\quad\text{for every $x\in\mathrm{supp}(V)$ and $t > 0$}.
  $$
  Then, by letting $t\to+\infty$, we deduce that $u\leq 0$ on $\mathrm{supp}(V)\neq \varnothing$, but this is clearly in contradiction with \eqref{eq:uzeropositive}. This completes the proof.
\end{proof}

\section{Optimality on model trees}\label{optimality}\setcounter{equation}{0}
We start by showing a general non-uniqueness criterium which holds for any graph $(G,\omega,\mu)$ such that \eqref{e7f} is fulfilled. We write 
$x\to \infty$ whenever $d(x, x_0)\to +\infty$, for some reference point $x_0\in G$.  

\subsection{A general non-uniqueness criterium}
Let us consider equation \eqref{problema}. We can prove the following result.
\begin{proposition}\label{prop71}
Let assumption \eqref{e7f} be in force. Moreover let $V\in\mathfrak{F}$, $V>0$ and $\hat R>0$. If there exists a supersolution to problem
\begin{equation}
\begin{aligned}\label{eq710}
\frac{1}{V}\Delta h=-1 \quad& \text{in}\,\,\, G\setminus B_{\hat R},\\
h>0\quad &\text{in}\,\,\,G\setminus B_{\hat R}, \\
 \lim_{x\to\infty}h(x)=0\,,&
\end{aligned}
\end{equation}
then there exist infinitely many bounded solutions $u$ of problem \eqref{problema}. In particular, for any $\gamma\in \R$, $\gamma>0$, there exists a solution $u$ to problem \eqref{problema} such that
$$
\lim_{x\to\infty} u(x)=\gamma.
$$
\end{proposition}

\begin{proof}
Let $\gamma\in\R$, $\gamma>0$. For any $j\in\mathbb{N}$, let us consider the following problem
\begin{equation}\label{eq711}
\begin{cases}
\Delta u-V(x)u&=0\quad \text{in}\,\,B_j\\
u&=\gamma \quad \text{in}\,\, G\setminus B_j.
\end{cases}
\end{equation}
Due to assumption \eqref{e7f}, existence and uniqueness of a solution $u_j$ to problem \eqref{eq711}, in the sense of Definition \ref{defsol}, for any $j\in\mathbb{N}$ is granted by Proposition \ref{prop:SolvDir}. We now claim that
\begin{equation}\label{eq712}
0\le u_j(x)\le \gamma \quad \text{for any}\,\,\,x\in G\,\,\,\text{and for any}\,\,\,j\in\mathbb{N}.
\end{equation}
In fact, since
$$
\Delta u_j-V(x)u_j=0\quad \text{in}\,\,B_j,
$$
and since $u_j=\gamma\ge0$ in $G\setminus B_j$, from Lemma \ref{lem:WMP}, we can infer that $u_j\ge0$ in $G$. On the other hand, let $\overline v(x):=\gamma$ for any $x\in G$. Then, since $V(x)>0$ for any $x\in G$, and since $\gamma>0$
\begin{equation}\label{eq713}
\begin{aligned}
\Delta \overline v(x)-V(x)\overline v(x)&=\frac{1}{\mu(x)}\sum_{y\in G}\omega(x,y)[\overline v(y)-\overline v(x)]-V(x)\overline v(x)\\
&=-\gamma\,V(x)\le 0\quad \text{for any}\,\,x\in G.
\end{aligned}
\end{equation}
For any $j\in \mathbb{N}$, let $w:=\overline v-u_j$. Due to \eqref{eq711} and \eqref{eq713},
$$
\Delta w-V(x) w\le 0\quad \text{in}\,\,G.
$$
Moreover, $w\ge0$ in $G\setminus B_{j}$. Hence, by Lemma \ref{lem:WMP}, $w\ge0$ in $G$ and, in particular,
$$u_j\leq  \gamma \quad \text{in}\,\, G\,.$$
Therefore, \eqref{eq712} follows. Furthermore, for any $j\in\mathbb{N}$, let $u_{j+1}$ be the solution to problem \eqref{eq711} in $B_{j+1}$. Thus, in particular, observe that
$$
\Delta u_{j+1}-V(x)u_{j+1}=0\quad \text{in}\,\,B_j,
$$
and, by \eqref{eq712}, $u_{j+1}(x)\le\gamma$ for any $x\in G\setminus B_j$. Let $v:=u_j-u_{j+1}$, then
$$
\begin{cases}
\Delta v-V(x)v=0\quad &\text{in}\,\,B_j,\\
v=\gamma-u_{j+1}\ge 0\quad &\text{in}\,\,G\setminus B_j.
\end{cases}
$$
Therefore, by Lemma \ref{lem:WMP}, $v\ge0$ in $G$ and, in particular, for any $j\in\mathbb{N}$,
\begin{equation}\label{eq714}
u_{j+1}\le u_j\quad \text{in}\,\,\,G.
\end{equation}
Hence, from \eqref{eq712} and \eqref{eq714}, we deduce that the sequence $\{u_j\}_{j\in\mathbb{N}}$ is \emph{decreasing and bounded} on $G$. Therefore, there exists $u\in\mathfrak{F}$ such that
$$u(x) = \lim_{j\to+\infty}u_j(x),\qquad x\in G.$$
Moreover, $u_j$ solves problem \eqref{eq711} for any $j\in\mathbb{N}$ and since the sum which defines the Laplacian $\Delta$ is finite, by letting $j\to+\infty$ we obtain that $u$ is a solution to problem \eqref{problema}.
Let $h$ be a supersolution to problem \eqref{eq710}. Then we define, for any $x\in G\setminus B_{\hat R}$, with $\hat R>0$ as in the assumptions,
\begin{equation}\label{eq715}
\underline w(x):=-C\,h(x)+\gamma.
\end{equation}
For any $j\in\mathbb{N}$ such that $B_j\supset B_{\hat R}$, we show that $\underline w$ is a subsolution to problem \eqref{eq711} in the sense of Definition \ref{defsol}.
Due to \eqref{eq710} and since $V>0$ in $G$, we have that
$$
\begin{aligned}
\Delta\underline w(x)-V(x) \underline w(x)&=\frac{C}{\mu(x)}\sum_{y\in G}\omega(x,y)[(-h(y)+\gamma)-(-h(x)+\gamma)]-V(x)(-Ch(x)+\gamma)\\
&=-C\Delta h(x)+V(x)C\,h(x)-V(x)\gamma\\
&\ge V(x)\left(C-\gamma\right)\\
&\ge 0\quad \text{for any}\,\,x\in B_j\setminus B_{\hat R}\,,
\end{aligned}
$$
provided that $C>0$ is big enough.
Moreover, $u_j\ge\underline w$ in $G\setminus B_j$ since $u_j=\gamma$ in $G\setminus B_j$. Therefore, by Lemma \ref{lem:WMP}, we get
\begin{equation}\label{eq716}
u_j(x)\ge \underline w(x)\quad \text{for any}\,\,x\in G\setminus B_{\hat R}.
\end{equation}
By combining together \eqref{eq712} and \eqref{eq716} we get
$$
-C\,h(x)+\gamma= \underline w(x)\le u_j(x)\le \gamma\quad \text{for any}\,\,x\in G\setminus B_{\hat R}.
$$
By letting $j\to\infty$ we get
$$
-C\,h(x)+\gamma\le u(x)\le \gamma\quad \text{for any}\,\,x\in G\setminus B_{\hat R},
$$
thus, in particular, due to \eqref{eq710},
$$
\lim_{x\to\infty} u(x)=\gamma.
$$
\end{proof}

\subsection{Model trees and a special supersolution}
In this subsection we consider a special kind of graphs, the so called \textit{model trees}, and we show that the uniqueness result in Theorem \ref{teo1} is sharp for this choice of graph. More precisely, we show that the choice $\alpha\in[0,1]$ in assumption \eqref{e12f} is optimal, indeed infinitely many bounded solutions exist whenever $\alpha>1$. Let us first define a model tree. Let $m(x)$ denote the number of edges which have $x$ as endpoint.
\begin{definition}\label{modeltree}
A graph $T$ will be called a \textit{model tree} if it contains a vertex $x_0$, known as the \textnormal{root} of the model, such that $m(x)$ is constant on spheres $S_r(x_0)=S_r$ of radius $r$ about $x_0$. Thus we have:
$$ \text{if}\,\,\, x\in S_r(x_0)=\{x\in T: \,d(x,x_0)=r\}, \text{ then}\quad m(x)=m(r).$$
\end{definition}
Here and hereafter, for any $x, x_0\in G$, $d(x, x_0)$ denotes the standard metric (on a graph), i.e. the number of vertices separating $x$ from $x_0$ along a path which connects $x_0$ to $x$.
  
Furthermore, we define the \textit{branching}, $b(r)$, at the distance $r$ from the root as the number of edges connecting each vertex in $S_r$ to a vertex in $S_{r+1}$; we also set $b(0)=m(x_0)$. Thus, for any $r>0$, we have that $m(r)=b(r)+1$. We denote a model tree with branching $b(r)$ by $T_{b(r)}$. We say that a model tree is \textit{homogeneous} if the branching is constant, i.e. $b(r)=b$ for every $r\ge0$, for some $b\in \mathbb N$.
\medskip

Let $\omega:G\times G\to\R$ be the edge weight of a graph $G$ as defined in \eqref{omega}. Whenever $\omega(x,y)\in\{0,1\}$ for all $x,y\in G$, we say that the graph has \textit{standard edge weight} and we denote it by $\omega_0(x,y)$.
In particular, for each $x\in G\setminus\{x_0\}$,
\begin{equation}\label{e3k}
\omega_0(x,y)=1   \text{ if and only if }  d(x,y)=1 \text{ and } [d(y, x_0)=d(x, x_0)\pm 1]\,.
\end{equation}
Observe that, for each $x\in G\setminus\{x_0\}$,
\begin{equation}\label{e4k}
\begin{aligned}
\operatorname{cardinality}\{y\in G\,: d(x,y)=1, d(y, x_0)=d(x, x_0)+1\}=b, \\
\operatorname{cardinality}\{y\in G\,: d(x,y)=1, d(y, x_0)=d(x, x_0)-1\}=1\,.
\end{aligned}
\end{equation}
Moreover, we define the weighted \textit{counting measure} as the measure $\mu$ for which there exists $c>0$ such that $\mu(x)=c$ for all $x\in G$, and we denote it by $\mu_c$. Moreover, let us denote the ball $B_r(x_0)$ of radius $r$ centered at the root $x_0$ simply by $B_r$.
\medskip

In what follows, we will deal with a \textit{model homogeneous tree with branching $b$, standard edge weight $\omega_0$ and weighted counting measure $\mu_c$}, that is the triplet
$$
(T_{b},\omega_0,\mu_c).
$$

Now, on $(T_{b},\omega_0,\mu_c)$, we explicitly construct a function $h$ satisfying all the properties required in Proposition \ref{prop71}. 

\begin{lemma}\label{barrier}
Let $(T_{b},\omega_0,\mu_c)$ be a graph as above with branching $b\geq 2$, and let $V\in\mathfrak{F},\,V> 0$
on $T_b$. Assume that,
for some $\alpha>1$ and $R_0\geq 2$, we have
\begin{equation}\label{e8k}
V(x) \leq C_0\, d^{-\alpha}(x, x_0)\quad \text{ for any }\; x\in G\setminus B_{R_0}\,.
\end{equation}
Then there exist $\hat C>0$, $\hat R>0$ and $\beta>0$ such that
\[h(x):=\hat C d^{-\beta}(x, x_0), \quad x\in G\setminus B_{\hat R}\]
is a supersolution to problem \eqref{eq710}.
\end{lemma}
\begin{proof} Let $x\in G \setminus B_{R_0}$ be fixed. Setting
$\mathbf{d}(\cdot) := d(\cdot,x_0)$, and observing that
$$\mathbf{d}(y)\geq\mathbf{d}(x)-1\geq 1\quad \forall\,\,y\in T_b,\,y\sim x$$
(recall that $R_0\geq 2$ and see \eqref{e3k}), by the Mean Value Theorem we have
\begin{equation}\label{e6k}
\begin{split}
\Delta h(x)& =\frac {\hat C}{\mu_c}\sum_{y\in G}\omega_0(x,y) [\mathbf{d}^{-\beta}(y) -\mathbf{d}^{-\beta}(x) ]
\\
& =-\beta \frac {\hat C}{\mu_c}\sum_{y\in G}\omega_0(x,y)\xi^{-\beta-1}[\mathbf{d}(y) -\mathbf{d}(x) ],
\end{split}
\end{equation}
for some $\xi$ in between $\mathbf{d}(y)$  and $\mathbf{d}(x)$.
Clearly,
\begin{equation}\label{e5k}
\begin{aligned}
\sum_{y\in G}\omega_0(x,y)\xi^{-\beta-1}[\mathbf{d}(y) -\mathbf{d}(x) ]& = \sum_{\mathbf{d}(y)>\mathbf{d}(x)}\omega_0(x,y)\xi^{-\beta-1}[\mathbf{d}(y) -\mathbf{d}(x) ]\\
&\qquad +\sum_{\mathbf{d}(y)<\mathbf{d}(x)}\omega_0(x,y)\xi^{-\beta-1}[\mathbf{d}(y) -\mathbf{d}(x) ]\\
& =: \mathrm{S}^+ + \mathrm{S}^-.
\end{aligned}
\end{equation}
Moreover,  by exploiting \eqref{e3k}-\eqref{e4k} (and since $\beta > 0$), we have
\begin{equation}\label{e7k}
\begin{aligned}
\mathrm{(i)}\,\,\mathrm{S}^+ &\geq
\sum_{\mathbf{d}(y)>\mathbf{d}(x)}\omega_0(x,y)\mathbf{d}^{-\beta-1}(y)[\mathbf{d}(y) -\mathbf{d}(x) ]
= b\, \mathbf{d}^{-\beta-1}(y)= - b [\mathbf{d}(x)+1]^{-\beta -1}; \\[0.05cm]
\mathrm{(ii)}\,\,\mathrm{S}^- &\geq \sum_{\mathbf{d}(y)<\mathbf{d}(x)}\omega_0(x,y)\mathbf{d}^{-\beta-1}(x)[\mathbf{d}(y) -\mathbf{d}(x) ] = \mathbf{d}^{-\beta -1}(x).
\end{aligned}
\end{equation}
By combining \eqref{e6k}, \eqref{e5k} and \eqref{e7k}, we then obtain
\begin{equation}\label{e2k}
\begin{aligned}
\Delta h(x)\leq \frac {\hat C \beta}{\mu_c }\left\{ - b [\mathbf{d}(x)+1]^{-\beta -1} + \mathbf{d}^{-\beta -1}(x)  \right\}.
\end{aligned}
\end{equation}
We finally observe that, given any $\e>0$, there exists $\hat R>R_0$ such that
\begin{equation}\label{e7kbis}
1<\frac{[\mathbf{d}(x)+1]^{\beta +1}}{\mathbf{d}^{\beta+1}(x) }<1+\e
\quad \text{ whenever } \mathbf{d}(x)>\hat R;
\end{equation}
thus, if $x\in G\setminus B_{\hat{R}}\subseteq G\setminus B_{R_0}$,
by \eqref{e2k}, \eqref{e7kbis} and \eqref{e8k} we obtain
\[\frac 1{V(x)}\Delta h(x)\leq  \frac {\hat C \beta }{\mu_c C_0}
\mathbf{d}^{\alpha}(x) \mathbf{d}^{-\beta -1}(x)\left(- \frac b{1+\e} +1\right)\leq - \frac {\hat C  \hat R^{\alpha-\beta-1}\beta }{2C_0 \mu_c} = -1,\]
provided that
$$\e=\frac{2b-3}{3}>0,\; 0<\beta\leq \alpha-1,\; \hat C=\frac{2C_0 \mu_c}{\hat R^{\alpha-\beta-1}\beta}.$$
Note that $\e>0$ since $b\geq 2$, while $\beta>0$ since $\alpha>1$. This completes the proof.
\end{proof}

From Proposition \ref{prop71} and Lemma \ref{barrier} we obtain the next
\begin{corollary}\label{nonuni}
  Let $(T_{b},\omega_0,\mu_c)$ be a graph as above with branching $b\geq 2$, and let $V\in\mathfrak{F},\,V> 0$
on $T_b$. Assume that condition \eqref{e8k} is fulfilled. 

Then, 
there exist infinitely many bounded solutions $u$ of problem \eqref{problema}. In particular, for any $\gamma\in \R$, $\gamma>0$, there exists a solution $u$ to problem \eqref{problema} such that
$$
\lim_{x\to\infty} u(x)=\gamma.
$$
\end{corollary}

\subsection{A counterexample}
On account of Corollary \ref{nonuni}, we can easily show that the requirement $\alpha\in [0,1]$
in assumption \eqref{e12f} \emph{cannot be dropped}, that is, this assumption
is \emph{optimal} if one restricts to a particular class of graphs.
\vspace*{0.1cm}

To illustrate this fact, let $(T_b,\omega_0,\mu_1)$
be a homogeneous model tree \emph{with branching $b \geq 2$} (and weighted
counting measure $\mu_c\equiv 1$), and let $x_0$ be the root of the model.
Given any number $\alpha> 1$, we then consider the function $V\in\mathfrak{F}$ defined as follows:
$$V(x) := (1+d(x,x_0))^{-\alpha}$$
(where $d$ is the usual distance on trees).
Clearly,
assumptions \eqref{e1f} and \eqref{e7f} are satisfied in this context;
moreover, $V:T_b\to\R$ is a \emph{strictly positive potential} on $T_b$, but the
last condition is assumption \eqref{e12f} is obviously violated (since $\alpha > 1$).

We now observe that, since $\alpha > 1$, the series \eqref{eq:hpStrada3} is convergent
\emph{for every choice of $\Lambda\in(0,1)$}: in fact, recalling that
$\mu_1\equiv 1$, we have the following computation
\begin{align*}
 & \sum_{x\in G\setminus B_{1}(x_0)}e^{- \Lambda d^\alpha(x,x_0)}\mu(x) =
 \sum_{n = 1}^{\infty}\Big(\sum_{\{x:\,d(x,x_0) = n\}}
 e^{- \Lambda d^\alpha(x,x_0)}\mu(x)\Big) \\
 & \qquad\qquad = \sum_{n = 1}^{\infty}
 e^{- \Lambda n^\alpha}\mu(\{x:\,d(x,x_0) = n\})
 = \sum_{n = 1}^{\infty}b^ne^{-\Lambda n^\alpha},
\end{align*}
 and the series $\sum_n b^ne^{-\Lambda n^\alpha}$ converges by the Ratio Test.
 \vspace*{0.1cm}
 
 On the other hand, since $b \geq 2$ and since
 $V$ is a strictly positive potential on $T_b$
 satisfying condition \eqref{e8k} (with $R_0 = 1$ and $C_0 = 1$),
 we are entitled to apply Corollary \ref{nonuni}, ensuring
  that there exist \emph{infinitely many} non-trivial bounded solutions
 to equation \eqref{problema}. Hence, the condition $\alpha\in [0,1]$ is optimal in this context.

\bigskip
\bigskip

\noindent{\bf Acknowledgement}
The second author is funded by the Deutsche Forschungsgemeinschaft (DFG, German Research Foundation) - SFB 1283/2 2021 - 317210226. All authors are member of the ``Gruppo Nazionale per l'Analisi Matematica, la Probabilit\'a e le loro Applicazioni'' (GNAMPA) of the ``Istituto Nazionale di Alta Matematica'' (INdAM, Italy).




\begin{thebibliography}{999}

\bibitem{BCG} M. Barlow, T. Coulhon, A. Grigor'yan, {\it Manifolds and graphs with slow heat kernel decay},  Invent. Math. {\bf 144} (2001), 609-649\,.

\bibitem{BP1} S. Biagi, F. Punzo, {\it A Liouville-type theorem for elliptic equations with singular coefficients in bounded domains}, Calc. Var. Part. Diff. Eq. DOI: 10.1007/s00526-022-02389-z, (to appear)\,.

\bibitem{BP2} S. Biagi, F. Punzo, {\it A Liouville theorem for elliptic equations in divergence form with a potential}, preprint (2022).



\bibitem{CGZ} T. Coulhon, A. Grigor'yan, F. Zucca, {\it The discrete integral maximum principle and its applications}, Tohoku J. Math. {\bf 57} (2005), 559-587.



\bibitem{Deo} N. Deo, "Graph Theory with Applications to Engineering and Computer Science", Dover Publications, New York (2016)


\bibitem{EM} M. Erbar, J. Maas, {\it Gradient flow structures for discrete porous medium equations}, Discr. Contin. Dyn. Syst. {\bf 34} (2014), 1355-1374.

\bibitem{Grig3} A.  Grigor'yan, {\it Bounded solutions of the Schr\"odinger equation on noncompact Rieamnnian manifolds}, J. Soviet Math. {\bf 51}, 2340-2349 (1990)\,.


\bibitem{Grig}  A. Grigor'yan, \emph{Analytic and geometric background of recurrence and non-explosion of the Brownian
motion on Riemannian manifolds}, Bull. Amer. Math. Soc. \bf 36 \rm (1999), 135--249.

\bibitem{GrigHK} A. Grigor'yan, "Heat Kernel and Analysis on Manifolds", AMS/IP Studies in Advanced Mathematics, 47. American Mathematical Society, Providence, RI; International Press, Boston, MA, 2009.

\bibitem{Grig2} A. Grigor'yan, "Introduction to Analysis on Graphs", AMS University Lecture Series {\bf 71} (2018)\,.


\bibitem{GLY1} A. Grigor'yan, Y. Lin, Y. Yang, {\it Kazdan-Warner equation on graph}, Calc. Var. Part. Diff. Eq. {\bf 55} (2016), 1-13.

\bibitem{GLY2} A. Grigor'yan, Y. Lin, Y. Yang, {\it  Yamabe type equations on graphs}, J. Diff. Eq. {\bf 261} (2016), 4924-943.

\bibitem{GLY3} A. Grigor'yan, Y. Lin, Y. Yang, {\it Existence of positive solutions to some nonlinear equations on locally finite graphs}, Sci. China Math. {\bf 60} (2017), 1311-1324.

\bibitem{GT} A. Grigor'yan, A. Telcs, {\it Sub-Gaussian estimated of heat kernels on infinite graphs}, Duke Math. J. {\bf 109(3)} (2001), 451--510.

\bibitem{GHY} Q. Gu, X. Huang, Y. Sun, {\it Superlinear elliptic inequalities on weighted graphs}, preprint (2022).


\bibitem{HL} B. Hua, Y. Lin, {\it Stochastic completeness for graphs with curvature dimension conditions}, Adv. Math. {\bf 306} (2017), 279-302\,.

\bibitem{HMu} B. Hua, D. Mugnolo, {\it Time regularity and long-time behavior of parabolic p-Laplace equations on infinite graphs}, J. Diff. Eq. {\bf 259} (2015), 6162-6190.

\bibitem{HW} B. Hua, L. Wang, {\it Dirichlet $p-$Laplacian eigenvalues and Cheeger constants on symmetric graphs}, Adv. Math. {\bf  364} (2020), 106997\,.

\bibitem{Huang} X. Huang, \emph{On uniqueness class for a heat equation on graphs}, J. Math. Anal. Appl. {\bf 393} (2012), 377--388.

\bibitem{HKMW} X. Huang, M. Keller, J. Masamune, R.K. Wojciechowski, {\it A note on self-adjoint extensions of the Laplacian on weighted graphs}, J. Funct. Anal. {\bf 265}, (2913), 1556-1578\,.

\bibitem{HKS} X. Huang, M. Keller, M. Schmidt, {\it On the uniqueness class, stochastic completeness and volume growth for graphs}, Trans. Amer. Math. Soc. {\bf 373} (2020), 8861-8884\,.

\bibitem{KS} M. Keller, M. Schwarz, {\it The Kazdan-Warner equation on canonically compactifiable graphs}, Calc. Var. Part. Diff. Eq. {\bf 57} (2018), 1-18.

\bibitem{KLW} M. Keller, D. Lenz, R.K. Wojciechowski,  "Graphs and Discrete Dirichlet Spaces", Springer (2021)\,.


\bibitem{Les} A. Lesne, {\it Complex networks: from graph theory to biology}, Lett. Math. Phys. {\bf 78}  (2006), 235-262\,.


\bibitem{LW2} Y. Lin, Y. Wu, {\it The existence and nonexistence of global solutions for a semilinear heat equation on graphs}, Calc. Var. Part. Diff. Eq. {\bf 56}, (2017), 1-22.

\bibitem{LY} S. Liu, Y. Yang, {\it Multiple solutions of Kazdan-Warner equation on graphs in the negative case}, Calc. Var. Part. Diff. Eq. {\bf 59} (2020), 1-15.


\bibitem{LHN} E. Lieberman, C. Hauert, M.A. Nowak, {\it Evolutionary dynamics on graphs}, Nature {\bf 433} (2005), 312-316.

\bibitem{MP} G. Meglioli, F. Punzo, \emph{Uniqueness for fractional parabolic and elliptic equations with drift}, Comm. Pure Applied Anal., to appear.

\bibitem{MP2} G. Meglioli, F. Punzo, \emph{Uniqueness in weighted $\ell^p$ spaces  for the Schr\"odinger equation on infinite graphs}, preprint (2022).


\bibitem{MR} G. Meglioli, A. Roncoroni, {\it Uniqueness in weighted Lebesgue spaces for an elliptic equation with drift on manifolds}, preprint (2022)

\bibitem{Mu} D. Mugnolo, {\it Parabolic theory of the discrete p-Laplace operator}, Nonlinear Anal. {\bf 87} (2013), 33-60.

\bibitem{Mu2} D. Mugnolo, "Semigroup Methods for Evolution Equations on Networks", Springer (2016)\,.




\bibitem{NP} C. Nobili, F. Punzo, {\it Uniqueness for degenerate parabolic equations in weighted $L^1$ spaces}, J. Evol. Eq. {\bf 22} (2022), 50\,.

\bibitem{PinaS} A. Pinamonti, G. Stefani, {\it Existence and uniqueness theorems for some semi-linear equations on locally finite graphs}, Proc. Amer. Math. Soc. {\bf 150} (2022), 4757-4770\,.


\bibitem{Pu1} F. Punzo, {\it Uniqueness for the heat equation in Riemannian manifolds}, J. Math. Anal. Appl. {\bf 424} (2015),
402-422.

\bibitem{Pu2} F. Punzo, {\it Integral conditions for uniqueness for solutions to degenerate parabolic equations}, J. Diff. Eq.
{\bf 267} (2019), 6555-6573.

\bibitem{PV} F. Punzo, E. Valdinoci, {\it Uniqueness in weighted Lebesgue spaces for a class of fractional parabolic and elliptic equations}, J. Diff. Eq. {\bf 258} (2015), 555-587.


\bibitem{SSV} A. Slavik, P. Stehlik, J. Volek, {\it Well-posedness and maximum principles for lattice reaction-diffusion equations}, Adv. Nonlinear Anal. {\bf 8} (2019), 303-322.


\bibitem{Wu} Y. Wu, {\it Blow-up for a semilinear heat equation with Fujita's critical exponent on locally finite graphs}, Rev. R. Acad. Cienc. Exactas Fis. Nat. Ser. A Mat. RACSAM {\bf 115} (2021), 1-16.

\end{thebibliography}
\end{document}